\documentclass[11pt,a4paper,twoside]{article}

\usepackage[T1]{fontenc}

\usepackage{amsmath,amssymb,latexsym,amsthm,varioref}
\usepackage{amsfonts,bbm,graphics,a4wide,rotating,graphicx,times,multirow,latexsym,fancybox}

\usepackage{graphicx,subfig,xcolor}
\labelformat{figure}{Figure~#1}

\renewcommand{\geq}{\ensuremath{\geqslant}}
\renewcommand{\leq}{\ensuremath{\leqslant}}
\newcommand{\N}{\ensuremath{\mathbb{N}}}
\newcommand{\Z}{\ensuremath{\mathbb{Z}}}
\newcommand{\R}{\ensuremath{\mathbb{R}}}
\newcommand{\1}{\ensuremath{\mathbf{1}}}
\renewcommand{\L}{\ensuremath{\mathbb{L}}}
\newcommand{\E}{\ensuremath{\mathbb{E}}}
\renewcommand{\P}{\ensuremath{\mathbb{P}}}
\def\var{\mathrm{Var}}
\newcommand{\f}{\ensuremath{\varphi_{\lambda}}}
\newcommand{\Nctot}{\ensuremath{N_{c,tot}}}

\numberwithin{equation}{section}
\allowdisplaybreaks

\theoremstyle{plain}  
\newtheorem{thm}{Theorem}
\labelformat{thm}{Theorem~#1}
\newtheorem{prop}{Proposition}
\labelformat{prop}{Proposition~#1}
\newtheorem{lem}{Lemma}[section]
\labelformat{lem}{Lemma~#1}
\theoremstyle{definition}  
\theoremstyle{remark}  

\title{\textbf{A model of Poissonian interactions \\
       and detection of dependence}}
\author{Laure \textsc{Sansonnet}\footnote{Corresponding author, Tel: (+33) 1 69 15 57 79, Fax: (+33) 1 69 15 60 34}
        \smallskip\\
        {\small{\emph{Laboratoire de Mathématiques, CNRS UMR 8628}}} \\
        {\small{\emph{Université Paris-Sud}}} \\
        {\small{\emph{15, rue Georges Clémenceau 91405 Orsay Cedex, France}}} \\
        {\small{E-mail: \texttt{laure.sansonnet@math.u-psud.fr}}}
        \smallskip\\
        and
        \smallskip\\
        Christine \textsc{Tuleau-Malot}
        \smallskip\\
        {\small{\emph{Univ.\,Nice Sophia Antipolis, CNRS, LJAD, UMR 7351}}} \\
        {\small{\emph{06100 Nice, France}}} \\
        {\small{E-mail: \texttt{malot@unice.fr}}}}
\date{December 3, 2013}

\begin{document}

\maketitle

\begin{center}
\begin{minipage}{0.75\textwidth}
\noindent \textbf{Abstract:}
This paper proposes a model of interactions between two point processes, ruled by a reproduction function $h$, which is considered as the intensity of a Poisson process. In particular, we focus on the context of neuroscience to detect possible interactions in the cerebral activity associated with two neurons. To provide a mathematical answer to this specific problem of neurobiologists, we address so the question of testing the nullity of the intensity $h$.
We construct a multiple testing procedure obtained by the aggregation of single tests based on a wavelet thresholding method. This test has good theoretical properties: it is possible to guarantee the level but also the power under some assumptions and its uniform separation rate over weak Besov bodies is adaptive minimax.
Then, some simulations are provided, showing the good practical behavior and the robustness of our testing procedure.

\medskip


\medskip

\end{minipage}
\end{center}

\section{Introduction}

In neuroscience, an important issue lies in a better understanding of the dynamics of cerebral activity in the cortex. In practice it is possible to measure, in vivo and for a specific task, the cerebral activity through the emission of action potentials by several neurons, and the specific interest of the neurobiologists is to understand how these action potentials appear. During a task, the recording of all arrival times of these action potentials (or spikes) on a neuron forms a spike train. From this point of view, the spike train can be modeled by a point process.

Several years ago it was thought that activities of different neurons during a task were independent (for example, see Barlow \cite{Bar72}); this explains why in the studies, the spike trains were usually modeled by independent Poisson processes. Today, thanks to technological advances in terms of recording brain activity, various studies show that this belief is false (for instance, see Gerstein \cite{Ger04} and Lestienne \cite{Les96}). Thus the recent studies consider neuronal assemblies instead of the separate neuronal activities. For example, activities of pairs of neurons, that have been recorded simultaneously, show that there exists a phenomenon called synchronization (for instance, see Grammont and Riehle \cite{GR99} and Grün \emph{et al.}\,\cite{GDA10}): the presence of a spike on one of the two spike trains can affect the apparition of a spike, with a delay, on the second spike train. From a biological point of view, such a phenomenon reflects a reality. Indeed, an action potential appears if the neuron is sufficiently excited. To obtain a sufficient excitation, two strategies exist: either the frequency of spikes received by a single neuron increases, or the receiving neuron receives less spikes but at the same time from different neurons. This second strategy is precisely the synchronization. From a biological point of view, it is less energy consuming and the reaction is faster. Therefore, the neurobiologists are interested in detecting the synchronization phenomenon. More generally, they want to detect whether or not neurons evolve independently of each other, a dependence being a hint of a functional connection during a task.

To mathematically answer this question, we need a model taking into account the possible interactions between two neurons. In neuroscience, a possible model is the Hawkes process (for example, see \cite{Haw71} for theoretical aspects and \cite{BM96,KRS10,PSCR11,PSCR12} for its introduction in neuroscience). The complete Hawkes process being, theoretically speaking, a very complicated model, we consider a modified version which is also realistic for the possible applications (in neuroscience, in genomics, \ldots) and for which it is possible to carry out computations. One possible model is the following one.
Let $N_p$ and $N_c$ be two point processes with respective intensity conditionally on the past
\begin{equation}
\tilde{\lambda}_p : t \longmapsto \mu_p \quad \mbox{and} \quad \tilde{\lambda}_c : t \longmapsto \mu_c + \int_{-\infty}^{t} h(t-u) \,dN_p(u), \label{def_condint}
\end{equation}
where $\mu_p>0$, $\mu_c>0$, $h:\R\rightarrow\R$ with $h(t)=0$ for $t\leq0$ and where $dN_p$ is the point measure associated with the process $N_p$. The parameters $\mu_p$ and $\mu_c$ describe the spontaneous part (in the context of neuroscience, the spontaneous apparition of spikes) and the function $h$ reflects the influence of $N_p$ on $N_c$. The function $\tilde{\lambda}_c$ which denotes the intensity conditionally on the past of $N_c$ specifically means that the probability that a new point appears on $N_c$ at time $t$ is the combination of the spontaneous part $\mu_c$ and the vote of each point of $N_p$ before $t$ through the function $h$.
Moreover, $N_p$ is a homogeneous Poisson process (for instance, see \cite{Kin93}) and $N_c$ is a special case of Hawkes process.
The biological problem which consists in knowing whether $N_p$ influences $N_c$ is equivalent to test the null hypothesis $\mathcal{H}_0$: "$h=0$" against the alternative $\mathcal{H}_1$: "$h\neq0$".

The above formulation of the intensity $\tilde{\lambda}_c$ is an integral form.
However it is possible conditionally on all the points of $N_p$ to have a vision in terms of descendants and no more in terms of intensity conditionally on the only past observations. Indeed, given $T$ a positive real number representing the time of record of the neuronal activity and given $n$ a fixed positive integer, conditionally on the event "the number of points of $N_p$ lying in $[0;T]$ is $n$", the points of the process $N_p$ obey the same law as a $n$-sample of uniform random variables on $[0;T]$, denoted $U_1,\ldots,U_n$ and named parents.
Thus, conditionally on $U_1,\ldots,U_n$, we can write $\tilde{\lambda}_c(t) = \mu_c + \sum_{i=1}^{n} h(t-U_i)$. This new expression of $\tilde{\lambda}_c$ can be interpreted as follows. Each $U_i$ gives birth independently to a Poisson process $N_c^i$ with intensity the function $h(t-U_i)$ with respect to the Lebesgue measure on $\R$, to which is added a homogeneous Poisson process $N_c^0$ with constant intensity $\mu_c$, representing the orphans. We consequently consider the aggregated process
\begin{equation}
N_c = \sum_{i=0}^{n} N_c^i \quad \mbox{whose intensity is given by the function} \quad \mu_c + \sum_{i=1}^{n} h(t-U_i) \label{def_aggrproc}
\end{equation}
and the points of the process $N_c$ are named children.
With this interpretation, the goal of the present paper is to test the "influence or not" of the parents on the children, via the reproduction function $h$.
This second writing contains many benefits. First, the assumption that the support of $h$ is included in $\R_+^*$ is not mandatory. With respect to the first formulation, this may appear like a minor difference, but in practice the impact is considerable. Indeed, if we refer to the context of neuroscience, assuming that the support of $h$ is in $\R_+$ means that one favors a sense of interactions, namely $N_p$ affects $N_c$. However in practice, we do not have this information a priori. Therefore, when the test does not reject $\mathcal{H}_0$, it means that $N_p$ does not seem to influence $N_c$, but it may be because in reality it is $N_c$ that affects $N_p$. We must be careful that the initially proposed model is not symmetric in terms of neurons and that a support in $\R_+$ does not really allow to answer the question of dependence. The causality is indeed represented by the fact that a child appears after its parent and therefore $h$ has to be supported in $\R_+$. Heuristically, a consequence is the following interpretation: if a parent has a child before its own birth, it may represent that the child is the parent and the parent the real child. Looking at both sides of the support (by considering $\R_+$ and also $\R_-$) makes the procedure in some sense adaptive to the causality of parent/child roles but it does not allow to symmetrize the test by inverting the parent/child roles. Indeed, in our model, one parent can have several children but a child has at most one parent.
Another advantage of this second writing is that it allows applications to other disciplines such as genomics where one studies for example the favored or avoided distances between patterns on a strand of DNA and where it is not always possible to know which pattern rules the other. More details about this application to genomics can be find in Sansonnet \cite{San12}, where the author proposes an estimation procedure of the function $h$, assumed to be well localized, based on wavelet thesholding methods, in a very similar model to the one studied here. The interested reader will find other estimation procedures of the function $h$ in this DNA context, by using a Hawkes' model in Gusto and Schbath \cite{GS05} and Reynaud-Bouret and Schbath \cite{RBS10}.

In this paper, given $T$ a positive real number representing the recording time and given $n$ a fixed positive integer, we consider a $n$-sample $(U_1,\ldots,U_n)$ of uniform random variables on $[0;T]$ representing the parents and we consider the model defined by (\ref{def_aggrproc}) for the children. For the simulation study, parents process $(U_i)_i$ is simulated according to a homogenous Poisson process of intensity $\mu_p$.
Since the null hypothesis $\mathcal{H}_0$: "$h=0$" means that conditionally on the total number of points of $N_c$, the points of the process $N_c$ are i.i.d.\,(independent and identically distributed) with uniform distribution, a first rather naive approach is to perform a Kolmogorov-Smirnov test (for instance, see \cite{Dar57}). But this test is not powerful, as illustrated in the section devoted to simulations.
The aim of this paper is then to build a more powerful and nonparametric test $\Phi_\alpha$ with values in $\{0,1\}$ of $\mathcal{H}_0$: "$h=0$" against the alternative $\mathcal{H}_1$: "$h\neq0$", rejecting $\mathcal{H}_0$ when $\Phi_\alpha=1$, with prescribed probabilities of first and second kind errors. The performance of the test $\Phi_\alpha$ is measured by its uniform separation rate (for example, see \cite{Bar02}).

In neuroscience, parametric methods exist to detect such dependence. For instance, the Unitary Event (UE) (see \cite{GDA10}) and the Multiple Tests based on a Gaussian Approximation of the UE (MTGAUE) (see \cite{TMRRBG12}) methods answer partially the problem by considering coincidences (see Section 4.4 for more details).
In the one-sample Poisson process model (that is to say $n=1$ and $\mu_c=0$ in our model), many papers deal with different problems of testing the simple hypothesis that an observed point process is a Poisson process with a known intensity.
We can cite for example the papers of Fazli and Kutoyants \cite{FK05} where the alternative is also a Poisson process with a known intensity, Fazli \cite{Faz07} where the alternatives are Poisson processes with one-sided parametric intensities or Dachian and Kutoyants \cite{DK06} where the alternatives are self-exciting point processes (namely, Hawkes processes).
In the nonparametric framework, Ingster and Kutoyants \cite{IK07} propose a goodness-of-fit test where the alternatives are Poisson processes with nonparametric intensities in a Sobolev $\mathcal{S}^{\delta}_2(R)$ or a Besov ball $\mathcal{B}^{\delta}_{2,q}(R)$ with $1 \leq q < \infty$ and known smoothness parameter $\delta$. They establish its uniform separation rate over a Sobolev or a Besov ball and show the adaptivity of their testing procedure in a minimax sense.

In some practical cases like the study of the expression of neuronal interactions or the study of favored or avoided distances between patterns on a strand of DNA, such smooth alternatives (Sobolev or Besov balls) cannot be considered. Indeed, the intensity of the Poisson process $N_c$ in these cases may burst at a particular position of special interest for the neuroscientist or the biologist. So we have to develop a testing procedure able to distinguish a constant function (or here a null function) from a function that has some small localized spikes. These features are not well captured by using classical Besov spaces.
Hence we focus in particular on alternatives based on sparsity rather than on alternatives based on smoothness. For this, we are interested in the computation of uniform separation rates over weak versions of Besov balls.
Such alternatives have already been considered. For instance, Fromont \emph{et al.}\,\cite{FLRB11} propose non-asymptotic and nonparametric tests of the homogeneity of a Poisson process that are adaptive over various Besov bodies simultaneously and in particular over weak Besov bodies.
Another example is Fromont \emph{et al.}\,\cite{FLRB12} which construct non-asymptotic and nonparametric multiple tests of the equality of the intensities of two independent Poisson processes, that are adaptive in the minimax sense over a large variety of classes of alternatives based on classical and weak Besov bodies in particular.

The test $\Phi_\alpha$ proposed in this paper consists in a multiple testing procedure obtained by aggregating several single tests based on a wavelet thresholding method as in Fromont \emph{et al.}\,\cite{FLRB11,FLRB12} (they also consider model selection and kernel estimation methods).
First, \ref{prop_1sterror} proves that the multiple test is an $\alpha$-level test and \ref{thm_multiproc} gives a condition on the alternative to ensure that our multiple test has a prescribed second kind error. This result reveals two regimes as in Sansonnet \cite{San12}. Indeed our model presents a double asymptotic through the number $n$ of parents and the recording time $T$ (namely, the length of the observations interval), which is not usual. Since $N_p$ is a homogeneous Poisson process with constant intensity $\mu_p$, the number $n$ of points of $N_p$ falling into $[0;T]$ is the realization of a Poisson random variable with parameter $\mu_pT$. As a consequence with very high probability, $T$ is proportional to $n$ and in this case, the uniform separation rates of the multiple test over weak Besov bodies are established by \ref{thm_rate}. Thus, our testing procedure is near adaptive in the minimax sense over a class of such alternatives.
The proofs of these results are essentially based on concentration inequalities (see \cite{Mas07}) and on exponential inequalities for $U$-statistics (see \cite{HRB03}).
Secondly, some simulations are carried out to validate our procedure from a practical point of view, which is compared with the classical Kolmogrov-Smirnov test, a test of homogeneity due to Fromont \emph{et al.}\,\cite{FLRB11} and a testing procedure proposed by Tuleau-Malot \emph{et al.}\,\cite{TMRRBG12}, which formalized a well-known procedure in neuroscience, namely the UE method (see Grün \emph{et al.}\,\cite{GDA10}).

The paper is organized as follows.
Section 2 deals with the description of our testing procedure.
Section 3 is devoted to the general results of the paper. The control of the probability of second kind error is ensured by \ref{thm_singleproc} for the single testing procedures and by \ref{thm_multiproc} for the multiple test. The uniform separation rates of the multiple test over weak Besov bodies are provided in \ref{thm_rate}.
Section 4 presents the simulation study.
The proofs of our main theoretical results are finally postponed in Section 6.

\section{Description of our testing procedure}

In the sequel, the support of $h$ is supposed to be compact and known. For instance, in neuroscience, there is a maximal time of synchronization (estimated to 40 ms) during a task according to the neuroscientists.
Without loss of generality, we suppose now that the support of $h$ is strictly included in $[-1;1]$ and that we observe the $U_i$'s (the parents) on $[0;T]$ and realizations of the process $N_c$ (the children) on $[-1;T+1]$.
In addition, we assume that $h$ belongs to $\L_1(\R)$ and $\L_\infty(\R)$ and consequently, we can consider the decomposition of $h$ on the Haar basis denoted by $\{\f,\lambda\in\Lambda\}$:
\[h = \sum_{\lambda\in\Lambda} \beta_\lambda \f \quad \mbox{with} \quad \beta_\lambda = \int_\R h(x)\f(x) \,dx,\]
where
\[\Lambda = \{\lambda = (j,k) : j \geq -1, k \in \Z\}\]
and for any $\lambda\in\Lambda$ and any $x\in\R$,
\[\f(x) = \left\{\begin{array}{cl} \phi(x-k) & \mbox{if $\lambda = (-1,k)$} \\ 2^{j/2} \psi(2^jx-k) & \mbox{if $\lambda = (j,k)$ with $j \geq 0$} \end{array}\right.,\]
with
\[\phi=\1_{[0;1]} \quad \mbox{and} \quad \psi=\1_{]\frac12;1]}-\1_{[0;\frac12]}.\]
The functions $\phi$ and $\psi$ are respectively the father and the mother wavelets.
Since the goal is to detect a signal, more precisely to detect if the function $h$ is identically null or not, the Haar basis is suitable in our context. Furthermore from a practical point of view, the use of the Haar basis yields fast algorithms, easy to implement. Nevertheless the theoretical results of the present paper can be generalized to a biorthogonal wavelet basis (see \cite{CDF92} for a definition of this particular basis) as in \cite{RBR10,RBRTM11,San12}.
We precise that we can easily extend our results to a function $h$ compactly supported in $[-A;A]$ for any $A>0$ by scaling the data by $\lceil A \rceil + 1$.

By considering this wavelet decomposition of $h$, the null hypothesis $\mathcal{H}_0$: "$h=0$" means that all the coefficients $\beta_\lambda$ are null and the alternative hypothesis $\mathcal{H}_1$: "$h\neq0$" means that there exists at least one non-zero coefficient. Since $h$ is strictly supported in $[-1;1]$, if one coefficient $\beta_{(-1,k)}$ is non-zero, then there exists at least one coefficient $\beta_{(j,k)}$ with $j \geq 0$ which is also non-zero. Therefore, we focus only on the coefficients $\beta_{(j,k)}$ with $j \geq 0$ and we introduce the following subset $\Gamma$ of $\Lambda$
\[\Gamma = \{\lambda = (j,k) \in \Lambda : j \geq 0, k \in \mathcal{K}_j\},\]
with $\mathcal{K}_j=\{k\in\Z : -2^j \leq k \leq 2^j-1\}$ ($\mathcal{K}_j$ is the set of integers $k$ such that the intersection of the support of $\f$ and $[-1;1]$ is not empty, with $\lambda=(j,k)$).

For every $\lambda$ in $\Gamma$, the coefficient $\beta_\lambda$ is estimated by
\[\hat{\beta}_\lambda = \frac{\mathcal{G}(\f)}{n}, \quad \mbox{with} \quad \mathcal{G}(\f) = \int_\R \sum_{i=1}^{n} \left[\f(x-U_i) - \frac{n-1}{n} \E_\pi(\f(x-U))\right] \,dN_c(x),\]
where $\pi$ denotes the uniform distribution on $[0;T]$ and $\E_\pi(f(U))$ the expectation of $f(U)$ where $U\sim\pi$ for any measurable function $f$.
These estimates, inspired by those proposed in \cite{San12} for a simpler model, namely with $\mu_c=0$, are unbiased:

\begin{prop}  \label{prop_unbest}
For all $\lambda=(j,k)$ in $\Gamma$, $\hat{\beta}_\lambda$ is an unbiased estimator of $\beta_\lambda$.
\end{prop}

The proof of \ref{prop_unbest} uses the fact that for all $\lambda$ in $\Gamma$, $\int_{-1}^{1} \f(t) \,dt = 0$, avoiding boundary effects (see Section 6.1).

In order to test the null hypothesis $\mathcal{H}_0$: "$h=0$" against $\mathcal{H}_1$: "$h\neq0$", namely "$\exists\lambda\in\Gamma, \beta_\lambda\neq0$", we first propose to test for all $\lambda\in\Gamma$, the null hypothesis $\mathcal{H}_0$ against the alternative $\mathcal{H}^\lambda_1$: "$\beta_\lambda\neq0$".
For each $\lambda\in\Gamma$, the associated simple test actually consists in testing "$\beta_\lambda=0$" against "$\beta_\lambda\neq0$" or more precisely, in testing the absence of variation of the function $h$ on a small interval.
Then in a second time, we will aggregate these simple tests to test the nullity of $h$ on its complete support.

\subsection{The single testing procedures}

Let us fix some $\alpha\in]0;1[$ and $\lambda\in\Gamma$.
We want to construct an $\alpha$-level test of the null hypothesis $\mathcal{H}_0$: "$h=0$" against $\mathcal{H}_1^\lambda$: "$\beta_\lambda\neq0$", from the observation of the parents $U_1,\ldots,U_n$ and the realization of the Poisson process $N_c$. We notice first that the null hypothesis entails in particular that $\beta_\lambda=0$.

We introduce the testing statistic $\hat{T}_\lambda$ defined by
\[\hat{T}_\lambda = |\hat{\beta}_\lambda|.\]
Our single test consists in rejecting the null hypothesis when $\hat{T}_\lambda$ is too large and more precisely, when
\[\hat{T}_\lambda > q_\lambda^{[U_1,\ldots,U_n;\Nctot]}(\alpha),\]
where $\Nctot$ is the (random) number of points of the process $N_c$ falling into $[-1;T+1]$ and for any $m\in\N^*$, $q_\lambda^{[U_1,\ldots,U_n;m]}(\alpha)$ is the $(1-\alpha)$-quantile conditionally on $U_1,\ldots,U_n$ of
\begin{equation}
\hat{T}_{\lambda,m}^0 = \frac{1}{n} \left| \sum_{k=1}^{m} \sum_{i=1}^{n} \left[\f(V^0_k-U_i) - \frac{n-1}{n} \E_\pi\big(\f(V^0_k-U)\big)\right] \right|,  \label{def_statH0}
\end{equation}
with $(V^0_1,\ldots,V^0_m)$ a $m$-sample with uniform distribution on $[-1;T+1]$ (namely a $m$-sample of the process $N_c$ under $\mathcal{H}_0$). We can easily prove that conditionally on $U_1,\ldots,U_n$ and $\Nctot=m$, $\hat{T}_\lambda$ and $\hat{T}_{\lambda,m}^0$ have exactly the same distribution under $\mathcal{H}_0$.
Thus, the corresponding test function is defined by
\begin{equation}
\Phi_{\lambda,\alpha}=\1_{\hat{T}_\lambda>q_\lambda^{[U_1,\ldots,U_n;\Nctot]}(\alpha)}.  \label{def_testfct}
\end{equation}

\subsection{The multiple testing procedure}

Previously, testing procedures have been built based on each single empirical coefficient $\hat{\beta}_\lambda$. We propose in this subsection to consider a collection of empirical coefficients instead of a single one, and to define a multiple testing procedure by aggregating the corresponding single tests.

Let $\{w_\lambda,\lambda\in\Gamma\}$ be a collection of positive numbers such that $\sum_{\lambda\in\Gamma}e^{-w_\lambda}\leq1$. This set allows us to put weights to empirical coefficients according to their index $\lambda=(j,k)\in\Gamma$.
Given $\alpha\in]0;1[$, we consider the test which rejects $\mathcal{H}_0$ when there exists at least one $\lambda$ in $\Gamma$ such that
\[\hat{T}_\lambda > q_\lambda^{[U_1,\ldots,U_n;\Nctot]}(u_\alpha^{[U_1,\ldots,U_n;\Nctot]}e^{-w_\lambda}),\]
where
\begin{equation}
\begin{split}
& u_\alpha^{[U_1,\ldots,U_n;\Nctot]} \\
& = \sup\left\{u>0 : \P\left(\max_{\lambda\in\Gamma} \left(\hat{T}_{\lambda,\Nctot}^0 - q_\lambda^{[U_1,\ldots,U_n;\Nctot]}(ue^{-w_\lambda})\right) > 0 \, \Big| \, U_1,\ldots,U_n;\Nctot \right) \leq \alpha \right\}.
\end{split}  \label{def_ualpha}
\end{equation}
The corresponding test function is defined by
\begin{equation}
\Phi_\alpha = \1_{\max_{\lambda\in\Gamma} \left(\hat{T}_\lambda-q_\lambda^{[U_1,\ldots,U_n;\Nctot]}(u_\alpha^{[U_1,\ldots,U_n;\Nctot]}e^{-w_\lambda})\right)>0}.  \label{def_multitestfct}
\end{equation}

We mention that, since the set $\Gamma$ is infinite countable, the number of tests to be performed is infinite and this is not a problem from a theoretical point of view. But in practice, we have to perform a finite number of single tests and so, we will fix a maximal resolution level $j_0$ and we will carry out the single tests $\Phi_{\lambda,\alpha}$ for $\lambda=(j,k)$ in $\Gamma$ with $j \leq j_0$. The role of $u_\alpha^{[U_1,\ldots,U_n;\Nctot]}$ is crucial in particular to guarantee the level of the multiple test and consequently, this quantity depends on the chosen maximal resolution level $j_0$ when we consider a finite number of single tests.


In the next section, we study the properties of the single tests $\Phi_{\lambda,\alpha}$ defined by (\ref{def_testfct}) and the multiple test $\Phi_\alpha$ defined by (\ref{def_multitestfct}), through their probabilities of first and second kind errors.

\section{Main theoretical results}

\subsection{Probability of first kind error}

We constructed our single and multiple tests in such a way that the first kind error, which measures the probability that the test wrongly rejects the null hypothesis, is less than $\alpha$.

\begin{prop}  \label{prop_1sterror}
Let $\alpha$ be a fixed level in $]0;1[$.
Then the single test $\Phi_{\lambda,\alpha}$ defined by (\ref{def_testfct}) for any $\lambda\in\Gamma$ and the multiple test $\Phi_\alpha$ defined by (\ref{def_multitestfct}) are of level $\alpha$.
Furthermore, $u_\alpha^{[U_1,\ldots,U_n;\Nctot]}$ defined by (\ref{def_ualpha}) satisfies $u_\alpha^{[U_1,\ldots,U_n;\Nctot]}\geq\alpha$.
\end{prop}

This result shows that the tests are exactly of level $\alpha$, which is required for a test from a non-asymptotic point of view (namely $n$ and $T$ are not required to tend to infinity).

\subsection{Probability of second kind error}

The second kind error, which measures the probability that the test does not wrongly reject the null hypothesis is not fixed by the testing procedure, unlike the first kind error. We have to control the probability of second kind error in such a way that it is close to 0, in order to obtain powerful tests.
The following theorem brings out a condition which guarantees that the single tests have a prescribed second kind error.

We denote by $\P_0$ the distribution of the aggregated process $N_c$ under $\mathcal{H}_0$, $\P_h$ the distribution of $N_c$ whose intensity conditionally on $U_1,\ldots,U_n$ is given by the function $\mu_c + \sum_{i=1}^{n} h(t-U_i)$ for any alternative $h$ and by $\E_h$ the corresponding expectation.
Since $h$ belongs to $\L_1(\R)$ and $\L_\infty(\R)$, we introduce $R_1$ and $R_\infty$ two positive real numbers such that $\|h\|_1 \leq R_1$ and $\|h\|_\infty \leq R_\infty$.

\begin{thm}  \label{thm_singleproc}
Let $\alpha$, $\beta$ be fixed levels in $]0;1[$.
Let $\zeta$ and $\kappa$ be positive constants depending on $\beta$, $\mu_c$, $R_1$ and $R_\infty$.
For all $\lambda\in\Gamma$, let $\Phi_{\lambda,\alpha}$ be the test function defined by (\ref{def_testfct}).
Assume that
\begin{equation}
\begin{split}
|\beta_\lambda| & > \sqrt{\frac{2\zeta}{\beta}\left(\frac{1}{n} + \frac{1}{T} + \frac{2^{-j}n}{T^2}\right)} \\
& \quad + \kappa \left\{\sqrt{\ln{(2/\alpha)}}\left(\frac{1}{\sqrt{n}} + \frac{1}{\sqrt{T}} + \frac{2^{-j/2}\sqrt{n}}{T}\right) + \ln{(2/\alpha)}\left(\frac{\sqrt{j}}{n} + \frac{j2^{j/2}}{n^{3/2}} + \frac{2^{-j/2}}{nT}\right)\right\},
\end{split}  \label{thm_singleproc:condition}
\end{equation}
for $\lambda=(j,k)$.
Then,
\[\P_h(\Phi_{\lambda,\alpha}=0) \leq \beta.\]
\end{thm}

Note that the quantity $\frac{1}{n} + \frac{1}{T} + \frac{2^{-j}n}{T^2}$ that appears under the square root of the first term of the right hand side of (\ref{thm_singleproc:condition}) is of the same order as the upper bound of the variance of the estimates $\hat{\beta}_\lambda$ (see Proposition 1 of \cite{San12}). Consequently, the right hand side of (\ref{thm_singleproc:condition}) can be viewed as a standard deviation term, since the other terms are not asymptotically larger than the first term if we assume that $2^j \leq n^2/(\ln{n})^2$, where asymptotic means $\min(n,T)\rightarrow+\infty$.

\ref{thm_singleproc} means that if the coefficient $\beta_\lambda$ is far enough from 0, then the probability of second kind error is controlled. This result gives a threshold for $\beta_\lambda$ from which our associated single testing procedure is able to detect a signal and shows that its power is larger than $1-\beta$.
Furthermore, if we consider the regime "$T$ proportional to $n$" in order to compare our result with known asymptotic rates of testing, Condition (\ref{thm_singleproc:condition}) can be easily obtained for instance if $\beta_\lambda^2 > C/n$ by assuming that $2^j \leq n^2/(\ln{n})^2$, with $C$ a positive constant.

Now we are interested in the power of the multiple testing procedure and the following theorem gives a condition on the alternative in order to ensure that our multiple test has a prescribed second kind error.

For an orthonormal basis $\{\f,\lambda\in L\}$ of a finite dimensional subspace $S_L$ of $\L_2(\R)$, we denote by $D_L$ the dimension of $S_L$ (namely the cardinal of $L$) and by $h_L$ the orthogonal projection of $h$ onto $S_L$.

\begin{thm}  \label{thm_multiproc}
Let $\alpha$, $\beta$ be fixed levels in $]0;1[$.
Let $\Phi_\alpha$ be the test function defined by (\ref{def_multitestfct}).
Assume that there exists at least one finite subset $L$ of $\Gamma$ such that
\begin{equation}
\|h_L\|_2^2 > \Big(C_1 D_L + C_2 \sum_{\lambda\in L}w_\lambda\Big) \left[\frac{1}{n}+\frac{n}{T^2}\right] + \Big(C_3 D_L + C_4 \sum_{\lambda\in L}w_\lambda + C_5 \sum_{\lambda\in L}w_\lambda^2\Big) \left[\frac{j_L}{n^2} + \frac{j_L^22^{j_L}}{n^{3}} + \frac{1}{n^2T^2}\right],  \label{thm_multiproc:condition}
\end{equation}
where $j_L=\max\{j\geq 0 : (j,k)\in L \ \mbox{with} \ k\in\mathcal{K}_j\}$ and $C_1$, $C_2$, $C_3$, $C_4$ and $C_5$ are positive constants depending on $\alpha$, $\beta$, $\mu_c$, $R_1$ and $R_\infty$.
Then,
\[\P_h(\Phi_\alpha=0) \leq \beta.\]
\end{thm}

This theorem means that if there exists one subspace $S_L$ of $\L_2(\R)$ such that $h_L$ (the orthogonal projection of $h$ onto $S_L$) lies outside a small ball around 0, then the probability of second kind error is controlled. This result gives a threshold for the energy of $h_L$ from which our multiple testing procedure is able to detect a signal and shows that its power is larger than $1-\beta$.
Furthermore, if we consider the regime "$T$ proportional to $n$" in order to compare our result with known asymptotic rates of testing, Condition (\ref{thm_multiproc:condition}) can be easily obtained for example if $\|h_L\|_2^2 > C\times\left(D_L+\sum_{\lambda\in L}w_\lambda+\sum_{\lambda\in L}w_\lambda^2\right)/n$ by assuming that $2^{j_L} \leq n^2/(\ln{n})^4$, with $C$ a positive constant.
Then, the separation rate between the null and the alternative hypotheses is of order $D_L/n$, and this is typical for testing procedures based on a thresholding approach (for instance, see \cite{FLRB11,FLRB12}). Usually, nested tests (namely based on model selection) achieve a faster rate of separation of order $\sqrt{D_L}/n$ (for example, see \cite{Bar02,BHL03}). But these latter tests are not adaptive over weak Besov bodies.
Consequently, the separation rate established by \ref{thm_multiproc} leads to sharp upper bounds for the uniform separation rates over such particular classes of alternatives and so, our multiple testing procedure will be proved to be adaptive over particular classes of alternatives, based on weak Besov bodies.

\subsection{Uniform separation rates}

Given some $\alpha,\beta\in]0;1[$, an $\alpha$-level test $\Phi_\alpha$ defined by (\ref{def_multitestfct}) has previously been built, with a probability of second kind error at most equals to $\beta$ if Condition (\ref{thm_multiproc:condition}) is satisfied. Then, given a class $\mathcal{S}_\delta$ of alternatives $h$, it is natural to measure the performance of the test via its uniform separation rate $\rho(\Phi_\alpha,\mathcal{S}_\delta,\beta)$ over $\mathcal{S}_\delta$ (see \cite{Bar02}) defined by
\[\rho(\Phi_\alpha,\mathcal{S}_\delta,\beta) = \inf \left\{\rho>0 : \sup_{h\in\mathcal{S}_\delta,\|h\|_2>\rho} \P_h(\Phi_\alpha=0) \leq \beta\right\}.\]
In order to compare our result with known asymptotic rates of testing, we consider the regime "$T$ proportional to $n$" in this subsection.

We introduce for $\delta>0$, $R>0$ the Besov body
\[\mathcal{B}^{\delta}_{2,\infty}(R) = \left\{f \in \L_2(\R) : f = \sum_{\lambda\in\Lambda} \beta_\lambda \f, \quad \forall j \geq 0, \sum_{k\in\mathcal{K}_j} \beta_{(j,k)}^2 \leq R^2 2^{-2j\delta}\right\}.\]
We also consider a weaker version of the above Besov bodies defined for $p>0$, $R'>0$ by
\[\mathcal{W}^*_p(R') = \left\{f \in \L_2(\R) : f = \sum_{\lambda\in\Lambda} \beta_\lambda \f, \quad \sup_{s>0} s^p \sum_{\lambda\in\Gamma} \1_{|\beta_\lambda|>s} \leq R'^p\right\}.\]
Whereas the spaces $\mathcal{B}^{\delta}_{2,\infty}(R)$ constitute an ideal class to measure the regularity of the possible alternatives $h$, the spaces $\mathcal{W}^*_p(R')$ constitute an ideal class to measure the sparsity of a wavelet decomposed signal $h$.
Indeed, if $f = \sum_{\lambda\in\Lambda} \beta_\lambda \f \in \mathcal{W}^*_p(R')$, then the associated sequence $\beta=(\beta_\lambda)_{\lambda\in\Gamma}$ satisfies $\sup_{\ell\in\N^*} \ell^{1/p}|\beta|_{(\ell)} < \infty$, where the sequence $(|\beta|_{(\ell)})_\ell$ is the non-increasing rearrangement of $\beta$: $|\beta|_{(1)} \geq |\beta|_{(2)} \geq \ldots \geq |\beta|_{(\ell)} \geq \ldots$. This condition gives a polynomial control of the decreasing rate of the sequence $(|\beta|_{(\ell)})_\ell$. The smaller $p$ is, the sparser is the signal.
There exists an embedding between Besov and weak Besov balls:
\[\mathcal{B}^{\delta}_{2,\infty}(R) \subset \mathcal{W}^*_{\frac{2}{1+2\delta}}(r),\]
where the radius $r$ of the weak Besov ball depends on $\delta$ and $R$ (more precisely, $r=4^\delta R/\sqrt{2^{2\delta}-1}$).
See \cite{KP00,Riv04a,Riv04b} for more details and for extensions in a more general setting.
So, we consider in this paper such alternatives based on the intersection of Besov and weak Besov bodies, namely sparse functions with a small regularity, see below.

To evaluate the uniform separation rates, we choose the following collection of weights $\{w_\lambda,\lambda\in\Gamma\}$ defined by
\begin{equation}
w_\lambda = 2\big(\ln{(j+1)}+\ln{(\pi/\sqrt{6})}\big) + \ln{|\mathcal{K}_j|},  \label{def_weight}
\end{equation}
for any $\lambda=(j,k)\in\Gamma$, where $|\mathcal{K}_j|$ is the cardinal of $\mathcal{K}_j$ (here, $2^{j+1}$). With this choice, the collection of weights satisfies the condition $\sum_{\lambda\in\Gamma}e^{-w_\lambda}\leq 1$.
The following theorem gives the uniform separation rates over $\mathcal{B}^{\delta}_{2,\infty}(R)\cap\mathcal{W}^*_{\frac{2}{1+2\gamma}}(R')$, where the parameter $\delta$ measures the regularity and the parameter $\gamma$ the sparsity.

\begin{thm}  \label{thm_rate}
Let $\alpha$, $\beta$ be fixed levels in $]0;1[$.
Assume that $T$ is proportional to $n$.
Let $\Phi_\alpha$ be the test function defined by (\ref{def_multitestfct}) with the weights $w_\lambda$'s defined by (\ref{def_weight}).
Then, for any $\delta>0$, $\gamma>0$, $R>0$, $R'>0$, if $2\delta>\gamma/(1+2\gamma)$
\[\rho(\Phi_\alpha,\mathcal{B}^{\delta}_{2,\infty}(R)\cap\mathcal{W}^*_{\frac{2}{1+2\gamma}}(R'),\beta) \leq C \left(\frac{\ln{n}}{n}\right)^{\frac{\gamma}{1+2\gamma}},\]
with $C$ a positive constant depending on $\delta$, $\gamma$, $R$, $R'$, $\alpha$, $\beta$, $\mu_c$, $R_1$ and $R_\infty$.
\end{thm}

If $\delta>\gamma$, then the set $\mathcal{B}^{\delta}_{2,\infty}(R)\cap\mathcal{W}^*_{\frac{2}{1+2\gamma}}(R')$ is reduced to $\mathcal{B}^{\delta}_{2,\infty}(R)$ (given the above embedding between Besov and weak Besov balls) that only measures the regularity. Since we are interested in sparse functions (with a small regularity), this is not the purpose here. Then we restrain our interpretation to the case $\gamma\geq\delta$.
Note that \ref{thm_rate} holds for instance with $\delta=1/4$ and for all $\gamma>0$. In this case, $\delta=1/4$ corresponds to the small regularity mentioned previously. Consequently, the main index $\gamma$, the sparsity index, governs the rates of convergence.

Considering the regime "$T$ proportional to $n$", uniform separation rates of the test $\Phi_\alpha$ given by \ref{thm_rate} match the minimax separation rates established by Theorem 1 of Fromont \emph{et al.}\,\cite{FLRB11}, if $2\delta>\gamma/(1+2\gamma)$ and also $\delta<\gamma/2$ and $\gamma>1/2$. Consequently, \ref{thm_rate} illustrates the optimality of our testing procedure in the minimax setting.
Furthermore, the upper bound of uniform separation rates of our test $\Phi_\alpha$ over $\mathcal{B}^{\delta}_{2,\infty}(R)\cap\mathcal{W}^*_{\frac{2}{1+2\gamma}}(R')$ has already been obtained, up to a logarithmic term, for a wavelet thresholding estimation method proposed by Sansonnet \cite{San12} in a very similar context and more precisely, this is equal to the minimax estimation rates of the maxisets of the thresholding estimation procedure (see \cite{KP00,RBR10,Riv04b} for more details). This means that it is at least as difficult to test as to estimate over such classes of alternatives. Note that on Sobolev or classical Besov spaces, testing rates are usually faster than estimation rates.

\section{Simulation study}  \label{simulation}

The scope of this section is to study our testing procedure from a practical point of view. Thus we consider different simulated data sets on which we apply our procedure and three other methods: the conditional Kolmogorov-Smirnov (\textbf{KS}) test, a test of homogeneity (\textbf{H}) developed by Fromont \emph{et al.}\,\cite{FLRB11} and a Gaussian Approximation of the Unitary Events (\textbf{GAUE}) method developed by Tuleau-Malot \emph{et al.}\,\cite{TMRRBG12}. Then, Section~\ref{influence} addresses the sensitivity to the maximal resolution level $j_0$.

The programs related on the implementation of our testing procedure have been coded in \texttt{Scilab 5.2} (Scilab Enterprises S.A.S, Orsay, France) and are available upon request. The other methods have been implemented with programs and softwares previously used by the initial authors.

\subsection{Description of the data}

We create different data sets that are to a certain extent a reflection of a neurobiological reality. We consider the spike trains of two neurons $N_p$ and $N_c$ which are modeled by two point processes with respective conditional intensity $\tilde{\lambda}_p$ and $\tilde{\lambda}_c$ defined by (\ref{def_condint}).

For real spike trains it is not reasonable to postulate the stationarity of $N_p$ and $N_c$, i.e.\,$\mu_p$ and $\mu_c$ are constant and considering the same function $h$ on the entire recording period $[0;T]$ (see Grün \emph{et al.}\,\cite{GDA10}). But this assumption is quite feasible on smaller time ranges (see Grammont and Riehle \cite{GR99} and Grün \cite{Gru96}). However, to date, we have no algorithmic and statistical tool to clearly identify the stationarity ranges.
Several methods (UE and MTGAUE, see \cite{TMRRBG12} for example) propose to perform many tests on different small windows of time and to use a multiple testing procedure (for instance, see Benjamini and Hochberg \cite{BH95}) to combine them. Hence those methods can solve, at least in practice, this stationarity problem. The aim of this simulation study is not to show how our testing procedure can be incorporated in a Benjamini and Hochberg's approach, which lies outside the scope of the present paper, but to discuss the advantage of our method on one small window of time. This explains the use of the simulated data described below.

\begin{figure}[ht]
\begin{center}
\begin{tabular}{cccc}
\textbf{A} & \includegraphics[scale=0.5,angle=0]{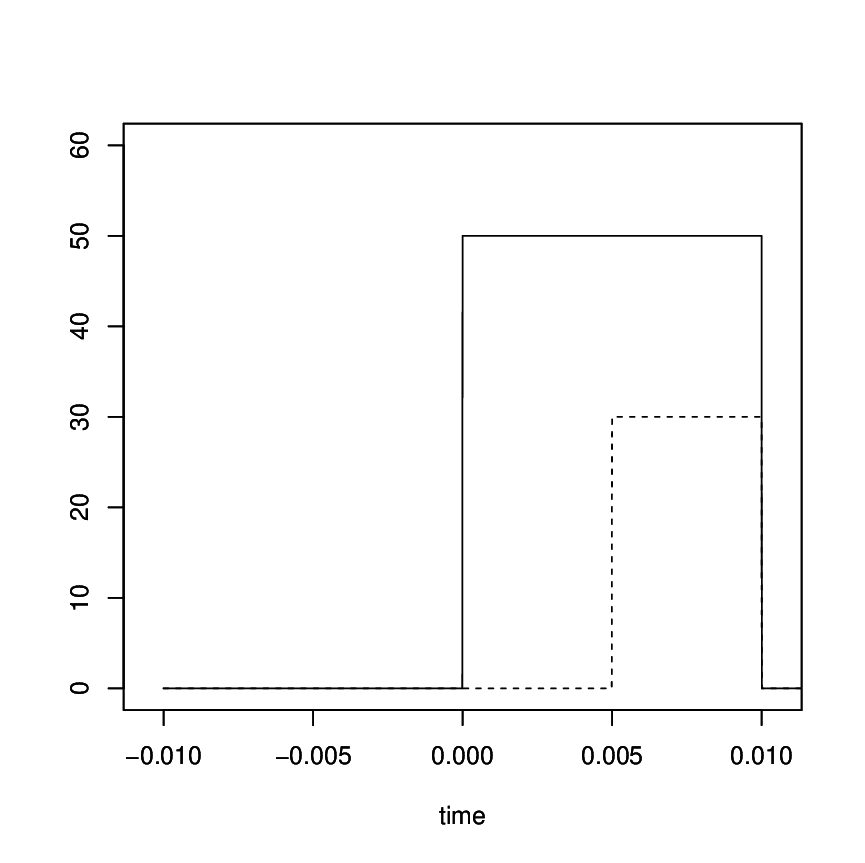} & \textbf{B} & \includegraphics[scale=0.5,angle=0]{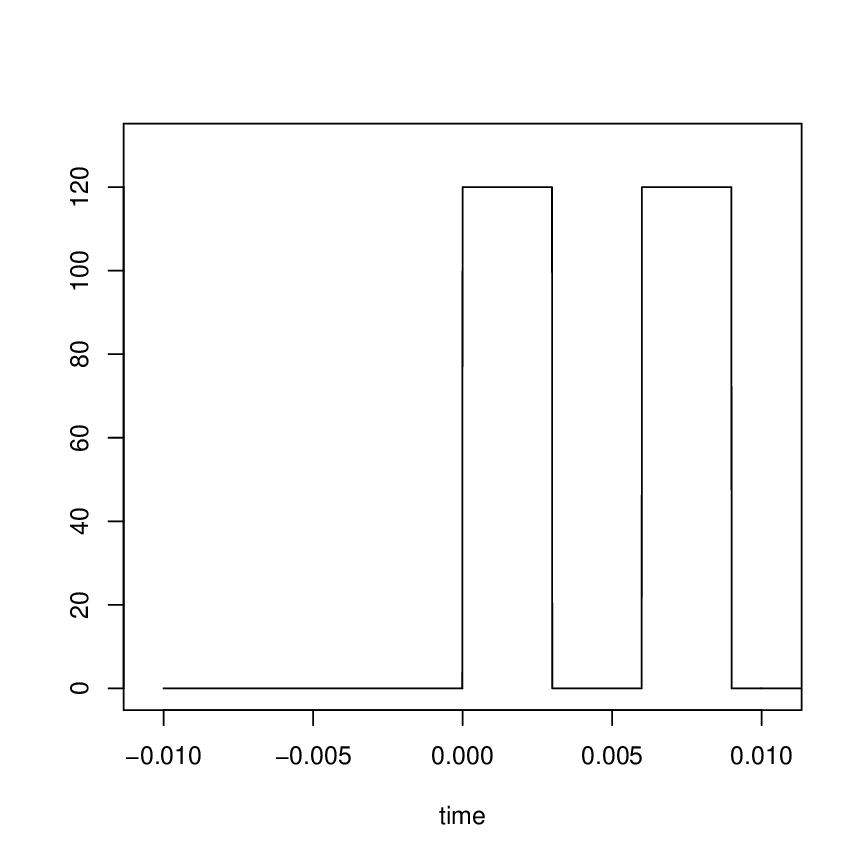}\\
\textbf{C} & \includegraphics[scale=0.5,angle=0]{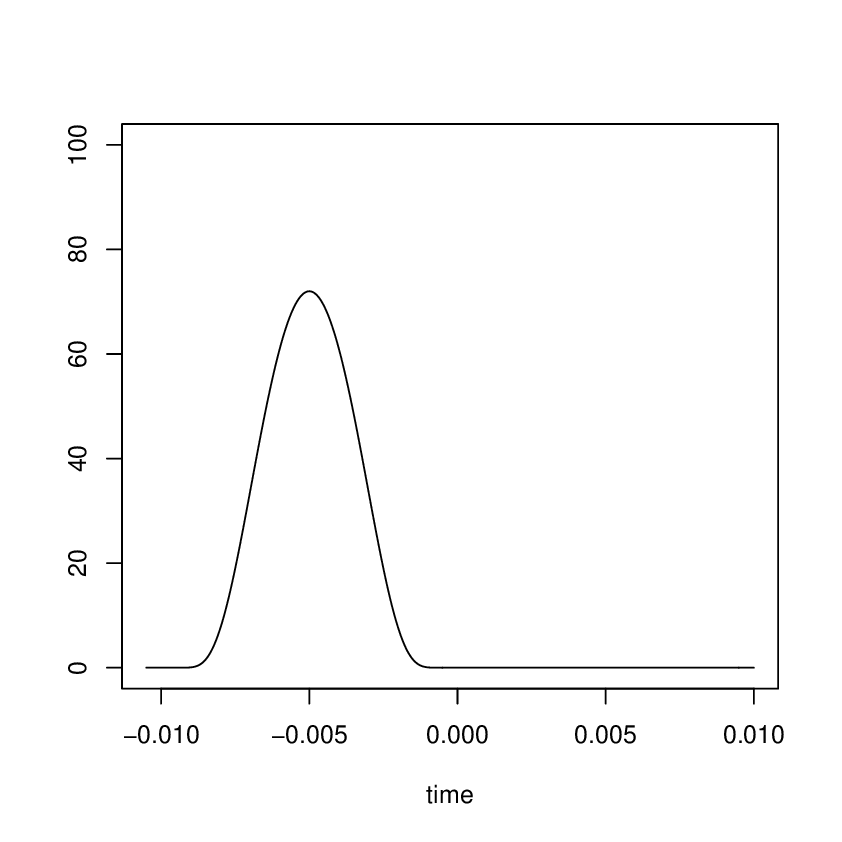} & \textbf{D} & \includegraphics[scale=0.5,angle=0]{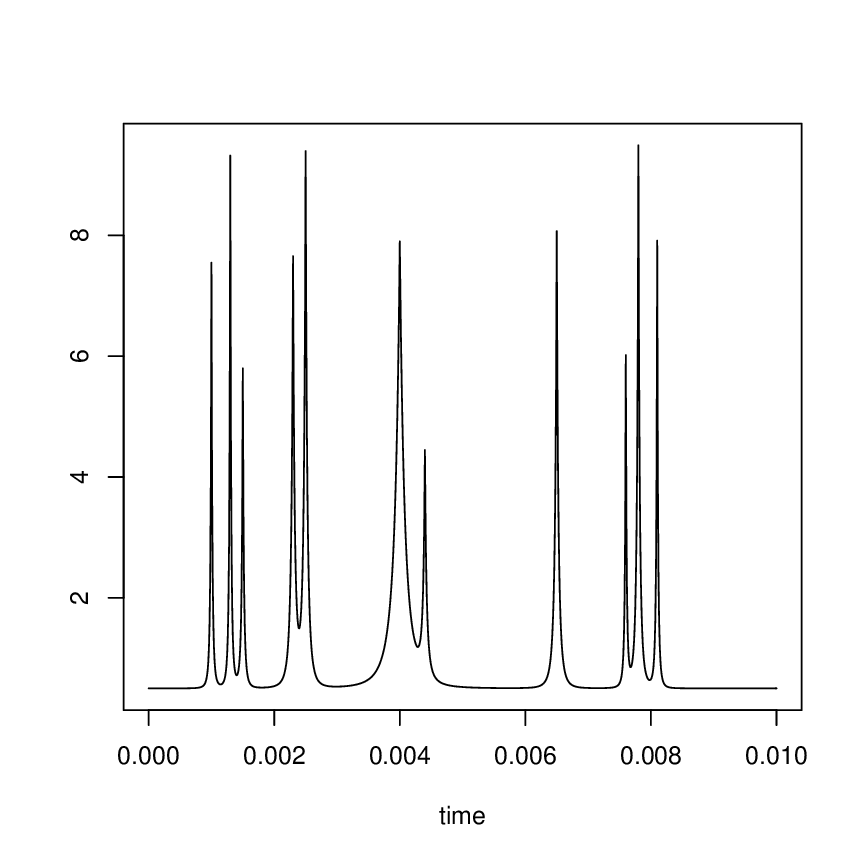}
\end{tabular}
\caption{{\small\emph{Graphs of alternatives: on Graph \textbf{A} the Echelon functions with $h_{50,0}$ in solid line and $h_{30,0.005}$ in dashed line, on Graph \textbf{B} the Crenel function $h_{Crenel}$, on Graph \textbf{C} the Bell function $h_{Bell}$ and on Graph \textbf{D} the Bumps function $h_{Bumps}$}}.}
  \label{fig_tf}
\end{center}
\end{figure}

We need therefore to simulate dependence between $N_p$ and $N_c$ on $[0;T]$, with $T=2\mbox{ s}$, and to take into account the major part of the neurobiological reality. So, we simulate processes $N_p$ and $N_c$ whose intensities are respectively given by
\begin{equation}
\tilde{\lambda}_p=50 \quad \mbox{and} \quad \tilde{\lambda}_c=50+\int_{-\infty}^{t} h(t-u) \,dN_p(u). \label{model_simu}
\end{equation}
At this stage, we can estimate the level of different procedures with $h\equiv0$ and in order to evaluate the powers of different procedures, several alternatives are tested. The first chosen alternative consists in intensities (Echelon functions) motivated by the context of neuroscience. Those intensities are defined by
\begin{center}
\begin{tabular}{lrcl}
Echelon functions & $h_{\theta,\nu}$ & $=$ & $\theta\1_{[\nu;0.01]}$,
\end{tabular}
\end{center}
with $\theta\in\{10,30,50,80\}$ and $\nu\in\{0,0.005\}$.
The parameter $\theta$ represents the influence strength of $N_p$ on $N_c$: the larger the parameter $\theta$ is, the higher the influence of $N_p$ on $N_c$ is. The parameter $\nu$ introduces a possible minimal delay in the synchronization, i.e.\,the synchronization of the neuronal activity occurs with a delay $\delta$ uniform on $[\nu;0.01]$.
To study the robustness of our procedure facing the other methods, we consider three other intensities (Crenel, Bell and Bumps) defined by
\begin{center}
\begin{tabular}{lrl}
Crenel function & $h_{Crenel}(x)=$ &  $120\left(\1_{[0;0.003]}(x)+\1_{[0.006;0.009]}(x)\right)$, \\
Bell function & $h_{Bell}(x)=$ & $72 \times \exp{\left(-4 \times \left(\frac{x+0.005}{0.005}\right)^2 \times (1-\left(\frac{x+0.005}{0.005}\right)^2)^{-1}\right)}\1_{[-1;0]}(x)$ , \\
Bumps function & $h_{Bumps}(x)=$ & $\frac12 \1_{[0;0.01]}(x)+\frac12\left(\sum_j g_j\left(1+\frac{|x-p_j|}{w_j}\right)^{-4}\right)\frac{\1_{[0;0.01]}(x)}{0.3}$,
\end{tabular}
\end{center}
where the vectors $g=(g_j)_j$, $p=(p_j)_j$ and $w=(w_j)_j$ are defined for example page 188 of \cite{FLRB11}.
These alternatives are represented in \ref{fig_tf}.

We mention that, with these different simulated data sets, we have in average 100 points for the process $N_p$ (the number of parents) and the average number of points of the process $N_c$ (the children) is given by Table \ref{nb_moy_enf} according to the different simulations.

\begin{table}[ht]
  \begin{center}
  \begin{tabular}{|c|c||c|c|}
  \hline
  Function $h$ & Average number of children & Function $h$ & Average number of children \\
  \hline
  $h_{10,0}$ & 111 & $h_{10,0.005}$ & 104 \\
  $h_{30,0}$ & 130 & $h_{30,0.005}$ & 115 \\
  $h_{50,0}$ & 150 & $h_{50,0.005}$ & 125 \\
  $h_{80,0}$ & 180 & $h_{80,0.005}$ & 140 \\
  $h\equiv 0$ & 100 & $h_{Crenel}$ & 131 \\
  $h_{Bell}$ & 130 & $h_{Bumps}$ & 151  \\
  \hline
  \end{tabular}
  \caption{{\small\emph{Average numbers of children according to the choice of the alternative}}.}
  \label{nb_moy_enf}
  \end{center}
\end{table}

\subsection{The Kolmogorov-Smirnov test and a test of homogeneity}

A first naive approach is to perform the classical Kolmogorov-Smirnov test (see Darling \cite{Dar57}) to convince us that this commonly used test is not reliable in this context. Indeed, even if the \textbf{KS} test is not a test of independence, the \textbf{KS} test may provide an answer to the problem. Since as said before, under $\mathcal{H}_0$ and conditionally on $U_1,\ldots,U_n$ and $\Nctot=m$, the observations of $N_c$ are i.i.d.\,with common law the uniform distribution on $[-1;T+1]$, looking for the adequation of $N_c$ with this law could be an idea to detect the rejection of $\mathcal{H}_0$. So, the use of the \textbf{KS} test is relevant.

In the same spirit, we can also compare our procedure to an adaptive test of homogeneity based on model selection, proposed by Fromont \emph{et al.}\,\cite{FLRB11} which has been shown to be more powerful than \textbf{KS} (see \cite{FLRB11}).
This one tests the null hypothesis "$\tilde{\lambda}_c$ is a constant function on its support (typically $[0;1]$)" against the alternative hypothesis "$\tilde{\lambda}_c$ is not a constant function". The \textbf{H} test consists in the aggregation of single tests as in our procedure, based on an estimation of the squared $\L_2$-distance between the target function $\tilde{\lambda}_c$ and the set of constant functions. For a sake of clarity, we briefly give a summary of this \textbf{H} test.
Let $\{X_l, l=1,\ldots,m\}$ be the points of the process $N_c$, $J \geq 1$ and $S_J$ the subspaces generated by the subsets $\{\phi,\varphi_{\lambda}, \lambda \in \Lambda_J \}$, with $\Lambda_J=\{(j,k), j \in \{0,\ldots,J-1\}, \ k \in \{0,\ldots, 2^j-1 \} \}$.
Focusing on one model $S_J$, they introduce
$T_J=\sum_{\lambda \in \Lambda_J} \tilde{T}_{\lambda}$,
where $\tilde{T}_{\lambda}=C\times\sum_{l\neq l'=1}^{m} \varphi_{\lambda}(X_l)\varphi_{\lambda}(X_{l'})$ with $C$ an absolute positive constant and then they consider the following test statistics $T_{\alpha}=\sup_{J\in\mathcal{J}}(T_J-q_J^m(u_{J,\alpha}^m))$, where $\mathcal{J}$ is a finite subset of $\N^*$, $q_J^m(u_{J,\alpha}^m)$ is the $(1-u_{J,\alpha}^m)$-quantile of the distribution of $T_J | \Nctot=m$ and $u_{J,\alpha}^m$ is defined as in \cite{FLRB11}. Finally, the corresponding test function is $\Phi_{\alpha}=\1_{T_{\alpha}>0}$.

\subsection{The GAUE method adapted to our context}

Before comparing the methods, we briefly return to the principle of the \textbf{GAUE} method.
The aim of the \textbf{GAUE} method is to detect the dependence on a single window $[0;T]$. This method is based on the coincidences with delay. More precisely for the couple of processes $(N_p,N_c)$, we compute the number of coincidences with delay $\delta$ on $[0;T]$, i.e.\,the variable $X_T = \int_{[0;T]^2} \1_{|x-y| \leq \delta} \,dN_{p}(x) \,dN_{c}(y)$, that represents the number of pairs $(x,y)$ in $N_p \times N_c$ such that $|x-y| \leq \delta$. This tuning parameter $\delta$ varies on a regular grid of $[0.001;0.04]$ with a step 0.001.
Let us define $\hat{\lambda}_p=N_p([0;T])/T$ and $\hat{\lambda}_c=N_c([0;T])/T$ where $N_p([0;T])$ and $N_c([0;T])$ denote respectively the number of spikes of $N_p$ and $N_c$ among $[0;T]$. The quantities $\hat{\lambda}_p$ and $\hat{\lambda}_c$ are estimators of $\tilde{\lambda}_p$ and $\tilde{\lambda}_c$.

We reject the null hypothesis $\mathcal{H}_0$: "$h=0$" when $X_T \geq \hat{m}_0 + \hat{\sigma}u_{1-\alpha/2}$, where $\hat{m}_0 = \hat{\lambda}_p \hat{\lambda}_c (2T\delta-\delta^2)$, $\hat{\sigma}^2 = \hat{\lambda}_p \hat{\lambda}_c (2T\delta-\delta^2) + \hat{\lambda}_p \hat{\lambda}_c
\left(\hat{\lambda}_p+\hat{\lambda}_c\right) \left(\frac23\delta^3-\frac1T\delta^4\right)$ and $u_{1-\alpha/2}$ is the $(1-\alpha/2)$-quantile of a standard normal. This threshold comes from the theory developed in \cite{TMRRBG12} and is adapted to our context. The quantity $\hat{m}_0$ is a plug-in estimator of the expectation of $X_{T}$ under $\mathcal{H}_0$ and $\hat{\sigma}^2$ is an estimator of the variance.
It can be shown that under the assumptions "$N_p$ and $N_c$ are Poisson processes" and "$N_p$ and $N_c$ are stationary", this test is asymptotically of level $\alpha$.
Further details about the meaning of those different estimators are given in \cite{TMRRBG12}.

The \textbf{GAUE} method was developed jointly with a neurophysiologist and it fits in line the UE method developed by Grün and coauthors (for example, see \cite{Gru96} and \cite{GDA10}), which is a commonly used method in neuroscience. One of its main disadvantage is that $\delta$ has to be chosen beforehand. Part of the aim of this work is to propose a more adaptive method.

\subsection{Our procedure in practice}

From a theoretical point of view, the support of the function $h$, denoted $[-A;A]$, should be strictly included in $[-1;1]$. Furthermore, a theoretical choice of the maximal resolution level $j_0$ is given by the condition: $2^{j_0}\leq n^2/(\log{n}^4)$. However, in practice, a trade-off between the choice of $j_0$ and the value of $A$ should be made. For instance, if $h=\1_{[0;A]}$ and if the order of magnitude of $A$ is $2^{-J}$ or $1-2^{-J}$, with $J>j_0+1$, our procedure does not allow to detect locally the jump of $h$ at $A$. To compensate this problem, we could increase the value of $j_0$. But, the choice of $j_0$ is restricted by the theoretical upper bound and especially, a greater $j_0$ leads to an increase of the computational time (due in particular to the evaluation of the quantiles which requires many iterations). Consequently, we propose to scale the data in order to have $A$ close to 1/2. Since the considered data sets have been built with a function $h$ supported by $[-0.01;0.01]$, the data are multiplied by 50 before being treated with our method.

Let us recall that our test rejects $\mathcal{H}_0$ when there exists at least one $\lambda=(j,k)$ in $\Gamma$ with $j \leq j_0$ such that
\[\hat{T}_\lambda > q_\lambda^{[U_1,\ldots,U_n;\Nctot]}(u_\alpha^{[U_1,\ldots,U_n;\Nctot]}e^{-w_\lambda}),\]
where $j_0\geq 1$ denotes the maximal resolution level, $u_\alpha^{[U_1,\ldots,U_n;\Nctot]}$ is defined by (\ref{def_ualpha}) and the $w_\lambda$'s are given by (\ref{def_weight}).
Hence, for each observation of the process $N_c$ whose number of points is denoted by $\Nctot=m$, given the points of $N_p$ denoted $U_1,\ldots,U_n$, we estimate $u_\alpha^{[U_1,\ldots,U_n;m]}$ and the quantiles $q_\lambda^{[U_1,\ldots,U_n;m]}$ by classical Monte Carlo methods based on the simulations of $B$ independent sequences $\{V^b, 1 \leq b \leq B\}$, where $V^b =(V_1^b,\ldots,V_m^b)$ is a $m$-sample of uniform variables on $[-1;T+1]$ (i.e.\,the law of $N_c$ under $\mathcal{H}_0$, conditionally on $U_1,\ldots,U_n$ and $\Nctot=m$).
We fix $B=20000$ in the sequel since for larger values of $B$, the gain in precision for the estimates of $u_\alpha^{[U_1,\ldots,U_n;m]}$ and $q_\lambda^{[U_1,\ldots,U_n;m]}$ becomes negligible.
We define for any  $\lambda=(j,k)$ in $\Gamma$ with $j \leq j_0$, for $1 \leq b \leq B$:
\[\hat{T}_{\lambda,m}^{0,b} = \frac{1}{n} \left|\sum_{k=1}^{m} \sum_{i=1}^{n} \left[\f(V^b_k-U_i) - \frac{n-1}{n} \E_\pi\big(\f(V^b_k-U)\big)\right] \right|.\]
We compute these $\hat{T}_{\lambda,m}^{0,b}$'s with a cascade algorithm (see Mallat \cite{Mal89}).

Half of the $m$-samples is used to estimate the quantiles by putting in ascending order the $\hat{T}_{\lambda,m}^{0,b} $'s for any $\lambda$.
The other half is used to approximate the conditional probabilities occurring in (\ref{def_ualpha}).
Then, $u_\alpha^{[U_1,\ldots,U_n;m]}$ is obtained by dichotomy, such that the estimated conditional probability occurring in (\ref{def_ualpha}) is less than $\alpha$, but as close as possible to $\alpha$.

For the comparison of our testing procedure to the three other methods, we have arbitrarily chosen $j_0 = 3$. With such a choice, our procedure considers 15 single tests $\Phi_{\lambda,\alpha}$ involving wavelets whose support length is respectively 0.125, 0.25, 0.5 and 1. This allows us to make detections at the positions $m\times2^{-3}$ ($m$ in $\{0,\ldots,7\}$) with a range of $2^{-3}$.
Due to the scaling of the data in our procedure, we need to divide the positions and the range of the possible detections by 50. Consequently, in the real time, the positions and the range become $m\times0.0025$ ($m$ in $\{0,\ldots,7\}$) and $0.0025$.

\subsection{Results}





We compare our testing procedure and the other methods on the different simulated data sets.
First, we focus on the empirical rate of the type I error which is an approximation of the level of the tests. Thus, we simulate 5000 independent realizations of (\ref{model_simu}) with $h\equiv0$, simulations on which we perform the present method and the other ones with level $\alpha=0.05$.
On those data, we evaluate the empirical rate of type I error. Those results are summarized in Table \ref{table-level}: all the testing methods seem to have a correct level in practice. This means that the number of wrong rejections of $\mathcal{H}_0$ is well controlled.


\begin{table}[ht]
  \begin{center}
  \begin{tabular}{|c|c|c|c|}
  \hline
  \textbf{our procedure} & \textbf{GAUE} & \textbf{H} & \textbf{KS} \\
  \hline
  0.047 & 0.0446/0.0510/0.0548 & 0.0638 & 0.051 \\
  \hline
  \end{tabular}
  \caption{{\small\emph{Empirical rate of type I error associated with our procedure and the other methods (\textbf{GAUE}, \textbf{H} and \textbf{KS}). The theoretical level is $\alpha=0.05$. Since the \textbf{GAUE} method depends on the tuning parameter $\delta$, the given value is the minimum/median/maximum of the empirical rate over all the $\delta$}}.}
  \label{table-level}
  \end{center}
\end{table}

Secondly, we want to see if the number of wrong rejections of $\mathcal{H}_1$ is also controlled. We consider the power of the tests which is the proportion of correct rejections of $\mathcal{H}_0$. To evaluate the power of the tests, we simulate 1000 independent realizations of (\ref{model_simu}) with different alternatives (Echelon, Crenel, Bell and Bumps functions). The results of the empirical power are given by Table \ref{table-power}.

\begin{table}[ht]
  \begin{center}
  \begin{tabular}{|l|c|c|c|c|}
  \hline
  Alternatives & \textbf{our procedure} & \textbf{GAUE} & \textbf{H} & \textbf{KS} \\
  \hline
  $h_{10,0}$ & 0.134 & 0.068/0.1085/0.168 & 0.062 & 0.040 \\
  $h_{10,0.005}$ &0.076 & 0.047/0.0575/0.077& 0.074 & 0.054 \\
  $h_{30,0}$& 0.656 & 0.154/0.3795/0.707 & 0.095 & 0.051 \\
  $h_{30,0.005}$& 0.341 & 0.050/0.1415/0.277 & 0.073 & 0.059\\
  $h_{50,0}$ & 0.939 & 0.278/0.6645/0.953 & 0.179 & 0.087 \\
  $h_{50,0.005}$  & 0.712 & 0.053/0.2825/0.589& 0.091 & 0.053 \\
  $h_{80,0}$ & 0.995 & 0.451/0.9160/0.998 & 0.362 & 0.113\\
  $h_{80,0.005}$ & 0.975 & 0.048/0.4900/0.879& 0.135 & 0.073  \\
  $h_{Crenel}$ & 0.949 & 0.255/0.437/0.993 & 0.112 & 0.069 \\
  $h_{Bell}$ & 0.672 & 0.046/0.3275/0.742 & 0.085 & 0.053 \\
  $h_{Bumps}$ & 0.948 & 0.139/0.701/0.967 & 0.159 & 0.082\\
  \hline
  \end{tabular}
  \caption{{\small\emph{Empirical power associated with our procedure and the other methods (\textbf{GAUE}, \textbf{H} and \textbf{KS}), evaluated for various alternatives. The theoretical level is $\alpha=0.05$. Since the \textbf{GAUE} method depends on the tuning parameter $\delta$, the given value is the minimum/median/maximum of the empirical rate over all the $\delta$}}.}
  \label{table-power}
  \end{center}
\end{table}

The power of the \textbf{KS} test is very low, as expected. The test of homogeneity \textbf{H} developed by Fromont et al.\,\cite{FLRB11} has a higher power, but this one remains smaller than the power of the two other methods. Thus, tests of homogeneity are not sufficient to detect dependence as expected.

Our procedure and the \textbf{GAUE} method are comparable in terms of power, even though the Echelon functions $h_{\theta,0}$ are particularly adapted to the \textbf{GAUE} method. However for the Echelon functions $h_{\theta,0.005}$, our method seems to have better performance since the power is higher. By considering the empirical power values of Table \ref{table-power}, it seems that both methods can be used to detect dependence.

Moreover, if both methods are comparable in terms of performance, it remains that the testing procedure proposed in this paper has an advantage over the \textbf{GAUE} method. In fact, our method is statistically adaptive. Indeed, the parameter $\delta$ which appears in the \textbf{GAUE} method is not calibrated in practice.
In our method, we aggregate the single tests over $(j,k)$. So on one hand, we do not need to specify this parameter but just an upper bound $j_0$, the maximal resolution level: the method through weights (\ref{def_weight}), adapts to this unspecified parameter $(j,k)$. But on the other hand, by looking at the single tests $\Phi_{\lambda,\alpha}$ that have supported the rejection, we are able to partially recover an important information for the practitioner: the position ($k2^{-j}$) and the range ($2^{-j}$) of the influence. In fact, by looking only at this single testing procedure, we get an upper value for $0.01$ and a lower value for $\nu$ on the range of delay $\delta$ of synchronization. To obtain more precise estimations of the support of $h$, we can consider an estimate of $h$, for example the one proposed by Sansonnet \cite{San12}.
The capacity of our method to get an information on $\nu$ is due to the fact that for a resolution level $j$ we consider different positions $k$. This is not possible with the \textbf{GAUE} method. This explains why the results on the Echelon functions $h_{\theta,0.005}$ are better with our method.

\subsection{Sensitivity to the maximal resolution level $j_0$} \label{influence}

For the comparison of our testing procedure to the other methods, we have chosen arbitrarily the maximal resolution level $j_0=3$. In this subsection, we propose to study the influence of the choice of this maximal resolution level $j_0$ on our testing procedure.

Since mentioned before, when we consider a finite number of single tests, $u_\alpha^{[U_1,\ldots,U_n;\Nctot]}$, defined by (\ref{def_ualpha}), depends on the chosen maximal resolution level $j_0$. The automatic calibration of $u_\alpha^{[U_1,\ldots,U_n;\Nctot]}$ during the practical procedure allows to guarantee a global level $\alpha$ for the multiple test as it is illustrated in Table \ref{table-nivj0}. We mention that the calibrated $u_\alpha^{[U_1,\ldots,U_n;\Nctot]}$ in practice satisfies \ref{prop_1sterror}: $u_\alpha^{[U_1,\ldots,U_n;\Nctot]}\geq\alpha$.

\begin{table}[ht]
  \begin{center}
  \begin{tabular}{|c|c|c|c|c|c|}
  \hline
  $j_0$ & 1 & 2 & 3 & 4 & 5 \\
  \hline
  Empirical rate of type I & 0.0508 & 0.0488 & 0.047 & 0.0474 & 0.0438 \\
  \hline
  \end{tabular}
  \caption{{\small\emph{Empirical rate of type I error associated with our procedure with different maximal resolution levels $j_0$. The theoretical level is $\alpha=0.05$}}.}
  \label{table-nivj0}
  \end{center}
\end{table}

\begin{figure}[ht]
\begin{center}
 \includegraphics[scale=0.75,angle=0]{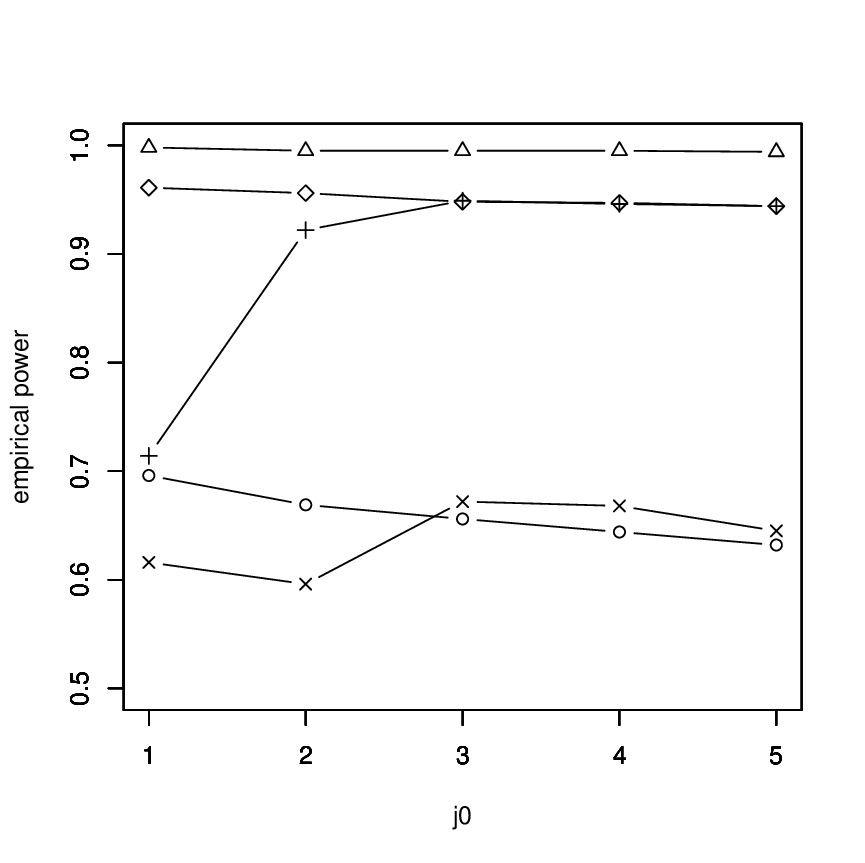}
 \caption{{\small\emph{Empirical power associated with our procedure according to $j_0$ for the alternatives $h_{30,0}$ in $\circ-\circ$, $h_{80,0}$ in $\vartriangle-\vartriangle$, $h_{Crenel}$ in $+-+$, $h_{Bell}$ in $\times-\times $ and $h_{Bumps}$ in $ \diamond-\diamond $. The theoretical level is $\alpha=0.05$}}.}
  \label{fig-j0}
\end{center}
\end{figure}

We are also interested in the influence of $j_0$ on the power of our test. \ref{fig-j0} displays the behavior of the power of our procedure according to the maximal resolution level $j_0$ for different alternatives.
We can first observe a stabilization of the power from $j_0=3$. Indeed, since the $w_\lambda$'s defined by (\ref{def_weight}) are not all identical and allocate different weights according to the index $\lambda=(j,k)$, weights decrease when the resolution level $j$ increases. Considering a higher maximal resolution level $j_0$ allocates a very small weight for the new tests of the procedure.
Furthermore, conforming to the real resolution level of the function which we want to test its nullity, we observe different behaviors for the first maximal resolution levels $j_0=1$ and $j_0=2$. For instance, the power of our procedure associated with the Crenel function is increasing with respect to $j_0$, whereas the power associated with the Echelon function $h_{30,0}$ is decreasing, but always with a kind of stabilization from $j_0=3$.

In light of this analysis of the influence of the maximal resolution level $j_0$ on our testing procedure, the choice of $j_0=3$ seems to be convenient, in order to obtain a suitable trade-off between power and computation time (we recall that the evaluation of the quantiles requires many iterations).

\section{Conclusion}

In our paper, we have investigated the influence of a point process on another one.
We have built a multiple testing procedure based on wavelet thresholding.
The main results of the paper have revealed the optimality of the procedure. Furthermore, our test is adaptive in the minimax sense over classes of alternatives essentially based on weak Besov bodies.
Then, from a practical point of view, our method answers several practical questions.
However, a number of challenges remain before applying our method on real data. To overcome the problem of stationarity, we could use a Benjamini and Hochberg's approach as for the GAUE method.
Finally, we could consider a more sophisticated model that takes into account the phenomenon of self-excitation (as for the complete Hawkes model). But this model raises serious difficulties from the theoretical point of view. This is an exciting challenge.

\section{Proofs}

All along the proofs, we introduce some positive constants denoted by $C(\xi,\ldots)$ meaning that they may depend on $\xi$, \ldots. They do not depend on $j$, $n$ and $T$ (which drive the asymptotic). Furthermore, the values of these constants may vary from line to line.

We recall that $\{\f,\lambda\in\Lambda\}$ is the Haar basis and consequently, we have:
\[\|\f\|_1 = 2^{-j/2}, \quad \|\f\|_2 = 1 \quad \mbox{and} \quad \|\f\|_\infty = 2^{j/2}.\]
In the case of a biorthogonal wavelet basis, $\|\f\|_1$, $\|\f\|_2$ and $\|\f\|_\infty$ are of the same order as above, up to a positive constant respectively depending on $\|\psi\|_1$, $\|\psi\|_2$ and $\|\psi\|_\infty$, where $\psi$ is the mother wavelet associated with the considered biorthogonal wavelet basis. Consequently, the same proofs potentially lead to the results on a biorthogonal wavelet basis as well as in \cite{San12} for the wavelet thresholding estimation.

\subsection{Proof of \ref{prop_unbest}}

We first notice that for any $\lambda$ in $\Gamma$, for any $u\in[0;T]$,
\begin{equation}
\int_{-1}^{T+1} \f(t-u) \,dt = 0.  \label{int_null}
\end{equation}
Let $\lambda\in\Gamma$ be fixed. By considering the aggregated process (\ref{def_aggrproc}), we can write
\begin{equation}
\mathcal{G}(\f) = G^0(\f) + G(\f),  \label{dec_G}
\end{equation}
with
\[G^0(\f) = \int_\R \sum_{i=1}^{n} \left[\f(x-U_i) - \frac{n-1}{n} \E_\pi(\f(x-U))\right] \,dN_c^0(x)\]
and
\[G(\f) = \int_\R \sum_{i=1}^{n} \left[\f(x-U_i) - \frac{n-1}{n} \E_\pi(\f(x-U))\right] \,\sum_{j=1}^{n} dN_c^j(x).\]
On the one hand, we notice that $G(\f)$ is the same quantity as the one defined by equation (2.2) of \cite{San12}. Thus, by applying the first part of Proposition 1 of \cite{San12}, we obtain
\[\E(G(\f)) = n \int_\R \f(x)h(x) \,dx.\]
On the other hand, we have
\[G^0(\f) = \int_\R \f(x-U_1) \,dN_c^0(x) + \sum_{i=2}^{n} \int_\R \left[\f(x-U_i)-\E_\pi(\f(x-U))\right] \,dN_c^0(x).\]
Thus,
\[\E(G^0(\f)|U_1,\ldots,U_n) = \int_{-1}^{T+1} \f(x-U_1) \mu_c \,dx + \sum_{i=2}^{n} \int_{-1}^{T+1} \left[\f(x-U_i)-\E_\pi(\f(x-U))\right] \mu_c \,dx\]
and by using (\ref{int_null}), we obtain
\[\E(G^0(\f))= \sum_{i=2}^{n} \int_{-1}^{T+1} \E\big[\f(x-U_i)-\E_\pi(\f(x-U))\big] \mu_c \,dx = 0.\]
Finally,
\[\E(\hat{\beta}_{\lambda}) = \E\left(\frac{\mathcal{G}(\f)}{n}\right) = \int_\R \f(x)h(x) \,dx = \beta_\lambda,\]
which proves \ref{prop_unbest}.

\subsection{Proof of \ref{prop_1sterror}}

Let $\alpha$ be a fixed level in $]0;1[$. Let $\lambda\in\Gamma$ be fixed.
First, the probability that the single test defined by (\ref{def_testfct}) wrongly detects a signal is
\[\P_0(\Phi_{\lambda,\alpha}=1) = \P_0\left(\hat{T}_\lambda>q_\lambda^{[U_1,\ldots,U_n;\Nctot]}(\alpha)\right).\]
Since conditionally on $U_1,\ldots,U_n$ and $\Nctot$, $\hat{T}_\lambda$ and $\hat{T}_{\lambda,\Nctot}^0$ have exactly the same distribution under $\mathcal{H}_0$, $q_\lambda^{[U_1,\ldots,U_n;\Nctot]}(\alpha)$ is also the $(1-\alpha)$-quantile of $\hat{T}_\lambda\big|U_1,\ldots,U_n;\Nctot$ under $\mathcal{H}_0$.
Thus,
\[\P_0(\Phi_{\lambda,\alpha}=1) \leq \alpha\]
and the level of the single test is $\alpha$.

Then, the probability that the multiple test defined by (\ref{def_multitestfct}) wrongly detects a signal is
\[\P_0(\Phi_\alpha=1) = \P_0\left(\max_{\lambda\in\Gamma} \left(\hat{T}_\lambda-q_\lambda^{[U_1,\ldots,U_n;\Nctot]}(u_\alpha^{[U_1,\ldots,U_n;\Nctot]}e^{-w_\lambda})\right)>0\right).\]
By definition (\ref{def_ualpha}) of $u_\alpha^{[U_1,\ldots,U_n;\Nctot]}$,
\[\P_0\left(\max_{\lambda\in\Gamma} \left(\hat{T}_\lambda-q_\lambda^{[U_1,\ldots,U_n;\Nctot]}(u_\alpha^{[U_1,\ldots,U_n;\Nctot]}e^{-w_\lambda})\right)>0 \Big| U_1,\ldots,U_n;\Nctot \right) \leq \alpha,\]
because conditionally on $U_1,\ldots,U_n$ and $\Nctot$, $\hat{T}_\lambda$ and $\hat{T}_{\lambda,\Nctot}^0$ have exactly the same distribution under $\mathcal{H}_0$.
By taking the expectation over $U_1,\ldots,U_n$ and $\Nctot$, we obtain that
\[\P_0(\Phi_\alpha=1) \leq \alpha\]
and the level of the multiple test is $\alpha$.

Furthermore, by Bonferroni's inequality we have
\begin{align*}
& \P\left(\max_{\lambda\in\Gamma} \left(\hat{T}_{\lambda,\Nctot}^0-q_\lambda^{[U_1,\ldots,U_n;\Nctot]}(\alpha e^{-w_\lambda})\right)>0 \Big| U_1,\ldots,U_n;\Nctot \right) \\
& \leq \sum_{\lambda\in\Gamma} \P\left(\hat{T}_{\lambda,\Nctot}^0-q_\lambda^{[U_1,\ldots,U_n;\Nctot]}(\alpha e^{-w_\lambda})>0 \Big| U_1,\ldots,U_n;\Nctot \right) \\
& \leq \sum_{\lambda\in\Gamma} \alpha e^{-w_\lambda} \\
& \leq \alpha
\end{align*}
and consequently $u_\alpha^{[U_1,\ldots,U_n;\Nctot]} \geq \alpha$ by definition (\ref{def_ualpha}) of $u_\alpha^{[U_1,\ldots,U_n;\Nctot]}$, which concludes the proof of \ref{prop_1sterror}.

\subsection{Proof of \ref{thm_singleproc}}

Let $\lambda\in\Gamma$ be fixed. Here we want to find a condition which will guarantee that
\[\P_h(\Phi_{\lambda,\alpha}=0)\leq\beta,\]
given $\beta\in]0;1[$.

Let us introduce $q^\alpha_{1-\beta/2}$ the $(1-\beta/2)$-quantile of the conditional quantile $q_\lambda^{[U_1,\ldots,U_n;\Nctot]}(\alpha)$.
Then for any $h$,
\begin{align*}
\P_h(\Phi_{\lambda,\alpha}=0) & = \P_h\left(\hat{T}_\lambda \leq q_\lambda^{[U_1,\ldots,U_n;\Nctot]}(\alpha) \, , \, q_\lambda^{[U_1,\ldots,U_n;\Nctot]}(\alpha) \leq q^\alpha_{1-\beta/2}\right) \\
& \quad + \P_h\left(\hat{T}_\lambda \leq q_\lambda^{[U_1,\ldots,U_n;\Nctot]}(\alpha) \, , \, q_\lambda^{[U_1,\ldots,U_n;\Nctot]}(\alpha) > q^\alpha_{1-\beta/2}\right) \\
&\leq \P_h(\hat{T}_\lambda \leq q^\alpha_{1-\beta/2}) + \beta/2
\end{align*}
and a condition which guarantees $\P_h(\hat{T}_\lambda \leq q^\alpha_{1-\beta/2}) \leq \beta/2$ will be enough to ensure that \[\P_h(\Phi_{\lambda,\alpha}=0) \leq \beta.\]
The following lemma gives such a condition.

\begin{lem}  \label{lem_2nderror}
Let $\alpha$, $\beta$ be fixed levels in $]0;1[$.
For any $\lambda=(j,k)\in\Gamma$, if
\begin{equation}
\E_h(\hat{T}_\lambda) > \sqrt{\frac{2\zeta Q_{j,n,T}}{\beta}} + q^\alpha_{1-\beta/2}  \label{lem_2nderror:condition}
\end{equation}
for a particular $\zeta$ which is a positive constant depending on $\mu_c$, $R_1$ and $R_\infty$,
where
\[Q_{j,n,T} = \frac{1}{n} + \frac{1}{T} + \frac{2^{-j}n}{T^2},\]
then
\[\P_h(\hat{T}_\lambda \leq q^\alpha_{1-\beta/2}) \leq \beta/2,\]
so that
\[\P_h(\Phi_{\lambda,\alpha}=0) \leq \beta.\]
\end{lem}

The proof of this lemma is postponed in Section 6.6.1.

In order to have an idea of the order of the right hand side of (\ref{lem_2nderror:condition}), we are now interested in the control of $q^\alpha_{1-\beta/2}$, the $(1-\beta/2)$-quantile of $q_\lambda^{[U_1,\ldots,U_n;\Nctot]}(\alpha)$. A sharp upper bound for $q^\alpha_{1-\beta/2}$ is given by the following lemma.

\begin{lem}  \label{lem_quant}
Let $\alpha$, $\beta$ be fixed levels in $]0;1[$.
For any $\lambda=(j,k)\in\Gamma$, there exists some positive constant $\kappa$ depending on $\beta$, $\mu_c$ and $R_1$ such that
\[q^{\alpha}_{1-\beta/2} \leq \kappa \left\{\sqrt{\ln{(2/\alpha)}} \left(\frac{1}{\sqrt{n}} + \frac{1}{\sqrt{T}} + \frac{2^{-j/2}\sqrt{n}}{T}\right) + \ln{(2/\alpha)} \left(\frac{\sqrt{j}}{n} + \frac{j2^{j/2}}{n^{3/2}} + \frac{2^{-j/2}}{nT}\right)\right\}.\]
\end{lem}

The proof of this lemma is postponed in Section 6.6.2.

Now, observe that if Condition (\ref{thm_singleproc:condition}) of \ref{thm_singleproc} is satisfied, namely
\[|\beta_\lambda| > \sqrt{\frac{2\zeta Q_{j,n,T}}{\beta}} + \kappa \left\{\sqrt{\ln{(2/\alpha)}}\left(\frac{1}{\sqrt{n}} + \frac{1}{\sqrt{T}} + \frac{2^{-j/2}\sqrt{n}}{T}\right) + \ln{(2/\alpha)}\left(\frac{\sqrt{j}}{n} + \frac{j2^{j/2}}{n^{3/2}} + \frac{2^{-j/2}}{nT}\right)\right\},\]
then by \ref{lem_quant},
\[|\beta_\lambda| > \sqrt{\frac{2\zeta Q_{j,n,T}}{\beta}} + q^{\alpha}_{1-\beta/2}.\]
We notice by Jensen's inequality that $|\beta_\lambda| = |\E_h(\hat{\beta}_\lambda)| \leq \E_h(|\hat{\beta}_\lambda|) = \E_h(\hat{T}_\lambda)$. Thus, Condition (\ref{lem_2nderror:condition}) of \ref{lem_2nderror} is satisfied and by \ref{lem_2nderror},
\[\P_h(\Phi_{\lambda,\alpha}=0) \leq \beta,\]
which concludes the proof of \ref{thm_singleproc}.

\subsection{Proof of \ref{thm_multiproc}}

Since $u_\alpha^{[U_1,\ldots,U_n;\Nctot]}\geq\alpha$ (see \ref{prop_1sterror}) and by setting $\alpha_\lambda=\alpha e^{-w_\lambda}$, we have
\begin{align*}
\P_h(\Phi_\alpha=0) & = \P_h\left(\forall\lambda\in\Gamma, \hat{T}_\lambda \leq q_\lambda^{[U_1,\ldots,U_n;\Nctot]}(u_\alpha^{[U_1,\ldots,U_n;\Nctot]}e^{-w_\lambda})\right) \\
& \leq \P_h\left(\forall\lambda\in\Gamma, \hat{T}_\lambda \leq q_{\lambda}^{[U_1,\ldots,U_n;\Nctot]}(\alpha_\lambda)\right) \\
& \leq \min_{\lambda\in\Gamma} \P_h\left(\hat{T}_\lambda \leq q_{\lambda}^{[U_1,\ldots,U_n;\Nctot]}(\alpha_\lambda)\right) \\
& \leq \min_{\lambda\in\Gamma} \P_h(\Phi_{\lambda,\alpha_\lambda}=0) \\
& \leq \beta,
\end{align*}
as soon as there exists $\lambda$ in $\Gamma$ such that $\P_h(\Phi_{\lambda,\alpha_\lambda}=0) \leq \beta$.

First, let us give the precise values of the constants that appear in Condition (\ref{thm_multiproc:condition}) of \ref{thm_multiproc}:
\[C_1 = 8\bigg(\frac{\zeta}{\beta}+3\kappa^2\ln{(2/\alpha)}\bigg), \ C_2 = 24\kappa^2, \quad C_3 = 12\kappa^2\ln^2{(2/\alpha)}, \ C_4 = 24\kappa^2\ln{(2/\alpha)} \ \mbox{and} \ C_5 = 12\kappa^2,\]
where $\zeta$ and $\kappa$ are the constants defined respectively by \ref{lem_2nderror} and \ref{lem_quant}.
We recall that $Q_{j,n,T}=\frac{1}{n}+\frac{1}{T}+\frac{2^{-j}n}{T^2}$ and we denote $R_{j,n,T}=\frac{j}{n^2} + \frac{j^22^{j}}{n^{3}} + \frac{2^{-j}}{n^2T^2}$.

Let us assume that there exists one finite subset $L$ of $\Gamma$ such that Condition (\ref{thm_multiproc:condition}) of \ref{thm_multiproc} is satisfied. Thus,
\begin{align*}
\|h_L\|_2^2 & > 8\left(\bigg(\frac{\zeta}{\beta}+3\kappa^2\ln{(2/\alpha)}\bigg) D_L + 3\kappa^2 \sum_{\lambda\in L}w_\lambda\right) \left[\frac{1}{n}+\frac{n}{T^2}\right] \\
& \quad + \left(12\kappa^2\ln^2{(2/\alpha)} D_L + 24\kappa^2\ln{(2/\alpha)} \sum_{\lambda\in L}w_\lambda + 12\kappa^2 \sum_{\lambda\in L}w_\lambda^2\right) \left[\frac{j_L}{n^2} + \frac{j_L^22^{j_L}}{n^{3}} + \frac{1}{n^2T^2}\right].
\end{align*}
Since $\ln{(2/\alpha)}+w_\lambda=\ln{(2/\alpha_\lambda)}$,
\[\sum_{\lambda\in L} \beta_\lambda^2 > \sum_{\lambda\in L} \left\{8\bigg(\frac{\zeta}{\beta}+3\kappa^2\ln{(2/\alpha_\lambda)}\bigg) \left[\frac{1}{n}+\frac{n}{T^2}\right] + 12\kappa^2\ln^2{(2/\alpha_\lambda)} \left[\frac{j_L}{n^2} + \frac{j_L^22^{j_L}}{n^{3}} + \frac{1}{n^2T^2}\right]\right\}\]
and it implies that there exists one coefficient $\lambda=(j,k)$ in $\Gamma$ such that
\[\beta_\lambda^2 > 8\bigg(\frac{\zeta}{\beta}+3\kappa^2\ln{(2/\alpha_\lambda)}\bigg) \left[\frac{1}{n}+\frac{n}{T^2}\right] + 12\kappa^2\ln^2{(2/\alpha_\lambda)} \left[\frac{j}{n^2} + \frac{j^22^{j}}{n^{3}} + \frac{1}{n^2T^2}\right].\]
Seeing that $Q_{j,n,T} \leq 2 \left[\frac{1}{n}+\frac{n}{T^2}\right]$ and $R_{j,n,T} \leq \left[\frac{j}{n^2} + \frac{j^22^{j}}{n^{3}} + \frac{1}{n^2T^2}\right]$, we have:
\[\beta_\lambda^2 > 4\frac{\zeta}{\beta} Q_{j,n,T} + 12\kappa^2\ln{(2/\alpha_\lambda)} Q_{j,n,T} + 12\kappa^2\ln^2{(2/\alpha_\lambda)} R_{j,n,T}.\]
Since $(\sqrt{a}+\sqrt{b}+\sqrt{c})^2\leq3(a+b+c)$ for all $a,b,c$ nonnegative real numbers,
\[\beta_\lambda^2 > 4\frac{\zeta}{\beta} Q_{j,n,T} + 4\kappa^2\ln{(2/\alpha_\lambda)} \left(\frac{1}{\sqrt{n}}+\frac{1}{\sqrt{T}}+\frac{2^{-j/2}\sqrt{n}}{T}\right)^2 + 4\kappa^2\ln^2{(2/\alpha_\lambda)} \left(\frac{\sqrt{j}}{n} + \frac{j2^{j/2}}{n^{3/2}} + \frac{2^{-j/2}}{nT}\right)^2\]
and then,
\[\beta_\lambda^2 > \left(\sqrt{\frac{2\zeta}{\beta} Q_{j,n,T}} + \kappa\bigg\{\sqrt{\ln{(2/\alpha_\lambda)}} \Big(\frac{1}{\sqrt{n}}+\frac{1}{\sqrt{T}}+\frac{2^{-j/2}\sqrt{n}}{T}\Big) + \ln{(2/\alpha_\lambda)} \Big(\frac{\sqrt{j}}{n} + \frac{j2^{j/2}}{n^{3/2}} + \frac{2^{-j/2}}{nT}\Big)\bigg\}\right)^2.\]
Finally, it is equivalent to
\[|\beta_\lambda| > \sqrt{\frac{2\zeta}{\beta} Q_{j,n,T}} + \kappa\left\{\sqrt{\ln{(2/\alpha_\lambda)}} \left(\frac{1}{\sqrt{n}}+\frac{1}{\sqrt{T}}+\frac{2^{-j/2}\sqrt{n}}{T}\right) + \ln{(2/\alpha_\lambda)} \left(\frac{\sqrt{j}}{n} + \frac{j2^{j/2}}{n^{3/2}} + \frac{2^{-j/2}}{nT}\right)\right\},\]
which is exactly Condition (\ref{thm_singleproc:condition}) of \ref{thm_singleproc} and we conclude the proof of \ref{thm_multiproc} by applying \ref{thm_singleproc}.

\subsection{Proof of \ref{thm_rate}}

With $T$ proportional to $n$, Condition (\ref{thm_multiproc:condition}) of \ref{thm_multiproc} is satisfied if there exists one finite subset $L$ of $\Gamma$ such that
\[\|h\|_2^2 \geq \|h-h_L\|_2^2 + C(\alpha,\beta,\mu_c,R_1,R_\infty) \bigg\{\Big(D_L + \sum_{\lambda\in L}w_\lambda\Big) \frac{1}{n} + \Big(D_L + \sum_{\lambda\in L}w_\lambda + \sum_{\lambda\in L}w_\lambda^2\Big) \left[\frac{j_L}{n^2} + \frac{j_L^22^{j_L}}{n^{3}}\right]\bigg\},\]
with $j_L=\max\{j\geq 0 : (j,k)\in L\}$, $\sum_{\lambda\in L}w_\lambda \leq C \times (j_L+1) D_L$ and $\sum_{\lambda\in L}w_\lambda^2 \leq C \times (j_L+1)^2 D_L$.
Consequently, Condition (\ref{thm_multiproc:condition}) is satisfied if there exists one finite subset $L$ of $\Gamma$ such that
\begin{equation}
\|h\|_2^2 \geq \|h-h_L\|_2^2 + C(\alpha,\beta,\mu_c,R_1,R_\infty) \frac{(j_L+1)}{n}D_L,  \label{pr-condition}
\end{equation}
with the maximal resolution level $j_L$ such that $2^{j_L} \leq n^2/(\ln{n})^4$.

Let $J\geq1$ that will be chosen later. We consider the following finite subset $\Gamma_J$ of $\Gamma$
\[\Gamma_J = \{\lambda=(j,k) \in \Gamma : 0 \leq j \leq J, k\in\mathcal{K}_j\}.\]
We introduce for all integer $D \leq |\Gamma_J|$ the subset $L$ of $\Gamma_J$ such that $\{\beta_\lambda, \lambda \in L\}$ is the set of the $D$ largest coefficients among $\{\beta_\lambda, \lambda \in \Gamma_J\}$.
We can notice that
\[\|h-h_L\|_2^2 = \|h-h_{\Gamma_J}\|_2^2 + \|h_{\Gamma_J}-h_L\|_2^2.\]
On the one hand, since $h$ belongs to $\mathcal{B}^{\delta}_{2,\infty}(R)$,
\[\|h-h_{\Gamma_J}\|_2^2 = \sum_{j>J}\sum_{k\in\mathcal{K}_{j}}\beta_{(j,k)}^2 \leq C(\delta)R^2 2^{-2J\delta}.\]
On the other hand, using equivalent definitions of weak Besov balls given by Lemma 2.2 of \cite{KP00} and using for instance page 211 of \cite{FLRB11}, we obtain:
\[\|h_{\Gamma_J}-h_L\|_2^2 \leq C(\gamma)R''^{2+4\gamma}D^{-2\gamma},\]
since $h$ belongs to $\mathcal{W}^*_{\frac{2}{1+2\gamma}}(R')$, with $R''$ an absolute positive constant depending eventually on $\gamma$ and $R'$.

Taking
\[J=\lfloor \log_2{(n^{\varepsilon})} \rfloor + 1\]
for some $0<\varepsilon<2$, we obtain that the right hand side of (\ref{pr-condition}) is upper bounded by
\[C(\delta,\gamma,R,R',\alpha,\beta,\mu_c,R_1,R_\infty) \left(n^{-2\varepsilon\delta}+D^{-2\gamma}+\frac{\varepsilon D \ln{n}}{n}\right).\]
Taking $D=\left\lfloor (n/\ln{n})^{1/(1+2\gamma)} \right\rfloor$ and $\varepsilon>\gamma/(\delta(1+2\gamma))$, we obtain that the right hand side of (\ref{pr-condition}) is upper bounded by
\[C(\delta,\gamma,R,R',\alpha,\beta,\mu_c,R_1,R_\infty) \left(\frac{n}{\ln{n}}\right)^{\frac{-2\gamma}{1+2\gamma}}\]
when $2\delta>\gamma/(1+2\gamma)$ and so,
\[\rho(\Phi_\alpha,\mathcal{B}^{\delta}_{2,\infty}(R)\cap\mathcal{W}^*_{\frac{2}{1+2\gamma}}(R'),\beta) \leq C(\delta,\gamma,R,R',\alpha,\beta,\mu_c,R_1,R_\infty) \left(\frac{n}{\ln{n}}\right)^{\frac{-\gamma}{1+2\gamma}},\]
which concludes the proof of \ref{thm_rate}.

\subsection{Proof of lemmas}

\subsubsection{Proof of \ref{lem_2nderror}}

Let $\lambda\in\Gamma$ be fixed. From Markov's inequality, we have that for any $x>0$,
\begin{equation}
\P_h\left(\left|\hat{T}_{\lambda}-\E_h(\hat{T}_{\lambda})\right| \geq x \right) \leq \frac{\var(\hat{T}_{\lambda})}{x^2}.  \label{pr-markov}
\end{equation}
Let us control $\var(\hat{T}_{\lambda})=\E_h(\hat{T}_{\lambda}^2)-\E_h^2(\hat{T}_{\lambda})$. We easily obtain by Jensen's inequality and by considering the decomposition (\ref{dec_G}) of $\mathcal{G}(\f)$:
\begin{align*}
\var(\hat{T}_{\lambda}) & \leq \var(\hat{\beta}_\lambda) \\
& \leq \frac{1}{n^2} \var(G^0(\f) + G(\f)) \\
& \leq \frac{2}{n^2} \big[\var(G^0(\f))+\var(G(\f))\big],
\end{align*}
with
\[\var(G(\f)) \leq C(R_1,R_\infty) \left\{n + \frac{n^2}{T} + \frac{2^{-j}n^3}{T^2}\right\},\]
by applying the second part of Proposition 1 of \cite{San12}.
It remains to compute $\var(G^0(\f))$. For this purpose, we apply the same methodology developed in Section 6.1.2 of \cite{San12}.
We have the following decomposition of $\var(G^0(\f))$ into two terms:
\begin{equation}
\var(G^0(\f))=\E(\var(G^0(\f)|U_1,\ldots,U_n)) + \var(\E(G^0(\f)|U_1,\ldots,U_n)).  \label{dec_varG0}
\end{equation}

We start by dealing with the first term of (\ref{dec_varG0}). We have
\begin{align*}
& \var(G^0(\f)|U_1,\ldots,U_n) \\
& = \int_{-1}^{T+1} \left(\sum_{i=1}^{n} \left[\f(x-U_i) - \frac{n-1}{n} \E_\pi(\f(x-U))\right]\right)^2 \mu_c \,dx \\
& = \mu_c \int_{-1}^{T+1} \left(\f(x-U_1) + \sum_{i=2}^{n} \left[\f(x-U_i)-\E_\pi(\f(x-U))\right]\right)^2 \,dx \\
& = \mu_c \int_{-1}^{T+1} \f^2(x-U_1) \,dx + 2\mu_c \int_{-1}^{T+1} \f(x-U_1) \sum_{i=2}^{n} \left[\f(x-U_i)-\E_\pi(\f(x-U))\right] \,dx \\
& \quad + \mu_c \int_{-1}^{T+1} \sum_{i=2}^{n} \sum_{k=2}^{n} \left[\f(x-U_i)-\E_\pi(\f(x-U))\right] \left[\f(x-U_k)-\E_\pi(\f(x-U))\right] \,dx.
\end{align*}
Since $\int_{-1}^{T+1}\f^2(x-U_1) \,dx = \|\f\|_2^2$,
\begin{align}
\E(\var(G^0(\f)|U_1,\ldots,U_n)) & = \mu_c \|\f\|_2^2 + \mu_c \int_{-1}^{T+1} \sum_{i=2}^{n} \E\left(\left[\f(x-U_i)-\E_\pi(\f(x-U))\right]^2\right) \,dx \nonumber\\
& = \mu_c \|\f\|_2^2 + (n-1)\mu_c \int_{-1}^{T+1} \var_\pi(\f(x-U)) \,dx \nonumber\\
& \leq \mu_c \|\f\|_2^2 + (n-1)\mu_c (T+2)\frac{\|\f\|_2^2}{T} \nonumber\\
& \leq C(\mu_c)n,  \label{pr-varG0-1}
\end{align}
by using (\ref{int_null}) and Lemma 6.1 of \cite{San12}.

Now, we deal with the second term of (\ref{dec_varG0}). We have
\begin{align*}
\E(G^0(\f)|U_1,\ldots,U_n) & = \int_{-1}^{T+1} \f(x-U_1) \mu_c \,dx + \sum_{i=2}^{n} \int_{-1}^{T+1} \left[\f(x-U_i)-\E_\pi(\f(x-U))\right] \mu_c \,dx \\
& = \mu_c \sum_{i=2}^{n} \int_{-1}^{T+1} \left[\f(x-U_i)-\E_\pi(\f(x-U))\right] \,dx,
\end{align*}
by using (\ref{int_null}).
Therefore,
\begin{align}
\var(\E(G^0(\f)|U_1,\ldots,U_n)) & = \mu_c^2 \var\left(\sum_{i=2}^{n} \int_{-1}^{T+1} \left[\f(x-U_i)-\E_\pi(\f(x-U))\right] \,dx\right) \nonumber\\
& = \mu_c^2 (n-1) \var\left(\int_{-1}^{T+1} \left[\f(x-U_1)-\E_\pi(\f(x-U))\right] \,dx\right) \nonumber\\
& \leq \mu_c^2 (n-1) \E\left[\left(\int_{-1}^{T+1} |\f(x-U_1)| \,dx\right)^2\right] \nonumber\\
& \leq \mu_c^2 (n-1) \|\f\|_1^2 \nonumber\\
& \leq C(\mu_c) 2^{-j}n.  \label{pr-varG0-2}
\end{align}

Finally, by combining inequalities (\ref{dec_varG0}), (\ref{pr-varG0-1}) and (\ref{pr-varG0-2}), we obtain:
\[\var(G^0(\f)) \leq C(\mu_c) n.\]

Thus,
\[\var(\hat{T}_{\lambda}) \leq \frac{C(\mu_c,R_1,R_\infty)}{n^2} \left\{n + \frac{n^2}{T} + \frac{2^{-j}n^3}{T^2}\right\} \leq \zeta Q_{j,n,T},\]
with
\[Q_{j,n,T} = \frac{1}{n} + \frac{1}{T} + \frac{2^{-j}n}{T^2}\]
and $\zeta$ a positive constant depending on $\mu_c$, $R_1$ and $R_\infty$.

Taking $x=\sqrt{2\zeta Q_{j,n,T}/\beta}$ in (\ref{pr-markov}) and using the previous inequality leads to
\[\P_h\left(\left|\hat{T}_{\lambda}-\E_h(\hat{T}_{\lambda})\right| \geq \sqrt{2\zeta Q_{j,n,T}/\beta} \right) \leq \frac{\beta}{2}.\]
Therefore, if $\E_h(\hat{T}_{\lambda}) > \sqrt{2\zeta Q_{j,n,T}/\beta} + q^\alpha_{1-\beta/2}$, then
\begin{align*}
\P_h(\hat{T}_{\lambda} \leq q^\alpha_{1-\beta/2}) & = \P_h\big(\hat{T}_{\lambda}-\E_h(\hat{T}_{\lambda}) \leq q^\alpha_{1-\beta/2}-\E_h(\hat{T}_{\lambda})\big) \\
& \leq \P_h\left(\left|\hat{T}_{\lambda}-\E_h(\hat{T}_{\lambda})\right| \geq \E_h(\hat{T}_{\lambda})-q^\alpha_{1-\beta/2}\right) \\
& \leq \P_h\left(\left|\hat{T}_{\lambda}-\E_h(\hat{T}_{\lambda})\right| \geq \sqrt{2\zeta Q_{j,n,T}/\beta}\right) \\
& \leq \beta/2
\end{align*}
and so
\[\P_h(\Phi_{\lambda,\alpha}=0) \leq \beta,\]
which concludes the proof of \ref{lem_2nderror}.

\subsubsection{Proof of \ref{lem_quant}}

We focus first on the control of the conditional quantile $q_\lambda^{[U_1,\ldots,U_n;\Nctot]}(\alpha)$. For all $m\in\N^*$, the $(1-\alpha)$-quantile $q_\lambda^{[U_1,\ldots,U_n;m]}(\alpha)$ is the smallest real number such that
\[\P\left(\hat{T}_{\lambda,m}^0 > q_\lambda^{[U_1,\ldots,U_n;m]}(\alpha) \bigg| U_1,\ldots,U_n;\Nctot=m\right) \leq \alpha,\]
where $\hat{T}_{\lambda,m}^0$ is defined by (\ref{def_statH0}).
Let $m\in\N^*$ be fixed. We write
\[\hat{T}_{\lambda,m}^0 = \frac{1}{n} \left|\sum_{k=1}^{m} S(\f)(V^0_k)\right|,\]
where $(V^0_1,\ldots,V^0_m)$ is a $m$-sample with uniform distribution on $[-1;T+1]$ and for any $v\in[-1;T+1]$,
\[S(\f)(v) = \sum_{i=1}^{n} \left[\f(v-U_i) - \frac{n-1}{n} \E_\pi(\f(v-U))\right].\]
Since $\E(\f(V-U)|U)=0$ for independent random variables $U$ and $V$ uniformly distributed on $[0;T]$ and $[-1;T+1]$ respectively, the $S(\f)(V^0_k)$'s are centered and independent conditionally on $U_1,\ldots,U_n$. Then we apply Bernstein's inequality (for instance, see Proposition 2.9 of \cite{Mas07}) to get that for all $\omega>0$, with probability larger than $1-2e^{-\omega}$,
\[\left|\sum_{k=1}^{m} S(\f)(V^0_k)\right| \leq \sqrt{2m\var(S(\f)(V^0_1)|U_1,\ldots,U_n)\omega} + \frac{\omega}{3} \sup_{v \in [-1;T+1]} \big|S(\f)(v)\big|.\]
Thus, with probability larger than $1-\alpha$,
\[\hat{T}_{\lambda,m}^0 \leq f(U_1,\ldots,U_n;m),\]
with
\begin{equation}
f(U_1,\ldots,U_n;m) = \frac{1}{n} \left\{\sqrt{2 m \ln{(2/\alpha)} V_S} + \frac{\ln{(2/\alpha)}}{3} B_S\right\},  \label{def_f}
\end{equation}
where
\[V_S=\var(S(\f)(V^0_1)|U_1,\ldots,U_n) \quad \mbox{and} \quad B_S=\sup_{v \in [-1;T+1]} |S(\f)(v)|.\]
Therefore we have $q_\lambda^{[U_1,\ldots,U_n;m]}(\alpha) \leq f(U_1,\ldots,U_n;m)$ by definition of the quantile $q_\lambda^{[U_1,\ldots,U_n;m]}(\alpha)$.

\medskip

Let us now provide a control in probability of $f(U_1,\ldots,U_n;m)$.
We control first $V_S$.
\begin{align}
V_S & = \var\left(\sum_{i=1}^{n} \f(V^0_1-U_i) - (n-1) \E_\pi\big(\f(V^0_1-U)\big) \Big| U_1,\ldots,U_n\right) \nonumber\\
& \leq \E\left[\Big(\sum_{i=1}^{n}\f(V^0_1-U_i)-(n-1)\E_\pi(\f(V^0_1-U))\Big)^2 \Big| U_1,\ldots,U_n\right] \nonumber\\
& \leq \frac{1}{T+2} \int_{v=-1}^{T+1} \left(\sum_{i=1}^{n}\f(v-U_i)-(n-1)\E_\pi(\f(v-U))\right)^2 \,dv \nonumber\\
& \leq \frac{2}{T+2} \int_{v=-1}^{T+1} \left(\sum_{1 \leq i,k \leq n}\f(v-U_i)\f(v-U_k) + (n-1)^2\E_\pi^2(\f(v-U))\right) \,dv \nonumber\\
& \leq \frac{2}{T+2} \Bigg\{\int_{v=-1}^{T+1} \sum_{i=1}^n\f^2(v-U_i) \,dv + \int_{v=-1}^{T+1} \sum_{1 \leq i \neq k \leq n}\f(v-U_i)\f(v-U_k) \,dv \nonumber\\
& \hspace{15em} + \frac{(n-1)^2}{T^2} \int_{v=-1}^{T+1} \left(\int_0^T |\f|(v-u) \,du\right)^2 \,dv\Bigg\} \nonumber\\
& \leq \frac{2}{T+2} \left\{n\|\f\|_2^2 + \int_{v=-1}^{T+1} \sum_{1 \leq i \neq k \leq n} \f(v-U_i)\f(v-U_k) \,dv + \frac{(n-1)^2}{T^2}(T+2)\|\f\|_1^2\right\} \nonumber\\
& \leq \frac{C}{T} \left\{n + \sum_{1 \leq i \neq k \leq n} \int_{v=-1}^{T+1} \f(v-U_i)\f(v-U_k) \,dv + \frac{2^{-j}n^2}{T}\right\},  \label{pr-majVS}
\end{align}
with $C$ an absolute positive constant.
We have a decomposition of the second term in a sum of degenerate $U$-statistics of order 0, 1 and 2. Indeed
\[\sum_{1 \leq i \neq k \leq n} \int_{v=-1}^{T+1} \f(v-U_i)\f(v-U_k) \,dv = W_0 + 2W_{1} + W_2,\]
with
\[W_{2}=\sum_{1 \leq i \neq k \leq n} \int_{v=-1}^{T+1} [\f(v-U_i)-\E_\pi(\f(v-U))][\f(v-U_k)-\E_\pi(\f(v-U))] \,dv,\]
\[W_{1}=\sum_{1 \leq i \neq k \leq n} \int_{v=-1}^{T+1} \f(v-U_i)\E_\pi(\f(v-U)) \,dv\]
and
\[W_{0}=-\sum_{1 \leq i \neq k \leq n} \int_{v=-1}^{T+1} \E_\pi^2(\f(v-U)) \,dv.\]

First we control $W_0$:
\begin{align}
|W_{0}| & \leq \frac{n(n-1)(T+2)}{T^2} \|\f\|_1^2 \nonumber\\
&\leq C \frac{2^{-j}n^2}{T},  \label{pr-W0}
\end{align}
with $C$ an absolute positive constant.
Next we deal with the control of $W_{1}$. We notice that
\[W_{1} = (n-1)\sum_{i=1}^{n} \int_{v=-1}^{T+1} \f(v-U_i)\E_\pi(\f(v-U)) \,dv\]
and consequently we have by using Lemma 6.3 of \cite{San12}
\begin{align}
|W_{1}| & \leq (n-1)\sum_{i=1}^{n} \int_{v=-1}^{T+1} |\f|(v-U_i) \,dv \frac{\|\f\|_1}{T} \nonumber\\
& \leq C \frac{2^{-j}n^2}{T},  \label{pr-W1}
\end{align}
with $C$ an absolute positive constant.

Now it remains to control $W_2$, with
\[W_2 = \sum_{1 \leq i < k \leq n} g(U_i,U_k),\]
where
\[g(U_i,U_k) = 2 \int_{v=-1}^{T+1} [\f(v-U_i)-\E_\pi(\f(v-U))][\f(v-U_k)-\E_\pi(\f(v-U))] \,dv.\]
One can apply Theorem 3.4 of \cite{HRB03} to $W_2$ and $-W_2$. It implies that there exist absolute positive constants $c_1$, $c_2$, $c_3$ and $c_4$ such that with probability larger than $1 - 2 \times 2.77 e^{-\omega}$,
\[|W_2| \leq c_1\mathcal{C}\sqrt{\omega} + c_2D\omega + c_3B\omega^{3/2} + c_4A\omega^2\]
for all $\omega > 0$, where
\begin{itemize}
  \item $A = \|g\|_{\infty} \leq 8 \|\f\|_1 \|\f\|_\infty \leq 8$;
  \item $\mathcal{C}^2 = \E(W_2^2)$ and we have
      \begin{align*}
      & \mathcal{C}^2 \\
      & = \sum_{1 \leq i < k \leq n} \E(g^2(U_i,U_k)) \\
      & \leq 4n(n-1) \E\left[\left(\int_{v=-1}^{T+1} [\f(v-U_1)-\E_\pi(\f(v-U))][\f(v-U_2)-\E_\pi(\f(v-U))] \,dv\right)^2\right] \\
      & \leq 4n(n-1) \E_{U,U'}\left[\left(\int_{v=-1}^{T+1} \big[\f^(v-U)-\E_\pi(\f(v-U))\big] \big[\f(v-U')-\E_\pi(\f(v-U'))\big] \,dv\right)^2\right] \\
      & \leq 4n^2 \E_{U,U'}\Bigg[\Bigg(\int_{v=-1}^{T+1} \f(v-U) \f(v-U') \,dv - \E_{U'}\left(\int_{v=-1}^{T+1} \f(v-U) \f(v-U') \,dv\right) \\
      & \hspace{6em} - \E_{U}\left(\int_{v=-1}^{T+1} \f(v-U) \f(v-U') \,dv\right) + \E_{U,U'}\left(\int_{v=-1}^{T+1} \f(v-U) \f(v-U') \,dv\right)\Bigg)^2\Bigg] \\
      & \leq Cn^2 \Bigg\{ \E_{U,U'}\left[\left(\int_{v=-1}^{T+1} |\f^U|(v) |\f^{U'}|(v) \,dv\right)^2\right] + \left[\E_{(U,U') \sim \pi \otimes \pi}\left(\int_{v=-1}^{T+1} |\f^U|(v) |\f^{U'}|(v) \,dv\right)\right]^2 \Bigg\},
      \end{align*}
      with $C$ an absolute positive constant. But,
      \begin{align*}
      & \E_{U,U'}\left[\left(\int_{v=-1}^{T+1} |\f|(v-U) |\f|(v-U') \,dv\right)^2\right] \\
      & \leq \E_{U,U'}\left(\int_{v=-1}^{T+1} |\f|^2(v-U) |\f|(v-U') \,dv \int_{v=-1}^{T+1} |\f|(v-U') \,dv\right) \\
      & = \E_{U,U'}\left(\int_{v=-1}^{T+1} |\f|^2(v-U) |\f|(v-U') \,dv\right) \|\f\|_1 \\
      & \leq \|\f\|_2^2 \frac{\|\f\|_1^2}{T}
      \end{align*}
      and
      \begin{align*}
      \E_{U,U'}\left(\int_{v=-1}^{T+1} |\f|(v-U) |\f|(v-U') \,dv\right) & = \int_{v=-1}^{T+1} \E_\pi(|\f|(v-U)) \E_\pi(|\f|(v-U')) \,dv \\
      & \leq (T+2)\frac{\|\f\|_1^2}{T^2},
      \end{align*}
      by using Lemma 6.3 of \cite{San12}. So,
      \[\mathcal{C}^2 \leq C n^2 \Bigg\{ \frac{2^{-j}}{T} + \frac{2^{-2j}}{T^2} \Bigg\},\]
      with $C$ an absolute positive constant;
  \item $\displaystyle D = \sup\left\{\E\left(\sum_{1 \leq k < i \leq n} g(U_i,U_k)a_i(U_i)b_k(U_k)\right) : \E\left(\sum_{i=2}^{n} a_i(U_i)^2\right) \leq 1, \E\left(\sum_{k=1}^{n-1} b_k(U_k)^2\right) \leq 1\right\}$. \\ But, with the conditions on the $a_i$'s and the $b_k$'s, we have:
      \begin{align*}
      & \E\left(\sum_{1 \leq k < i \leq n} g(U_i,U_k)a_i(U_i)b_k(U_k)\right) \\
      & = 2 \E\left(\sum_{i=2}^{n} \sum_{k=1}^{i-1} \int_{v=-1}^{T+1} [\f(v-U_i)-\E_\pi(\f(v-U))][\f(v-U_k)-\E_\pi(\f(v-U))] \,dv a_i(U_i)b_k(U_k)\right) \\
      & \leq 2 \int_{v=-1}^{T+1} \E\bigg(\sum_{i=2}^{n} \big|\f(v-U_i)-\E_\pi(\f^U(v))\big| |a_i(U_i)| \bigg) \E\bigg(\sum_{k=1}^{n-1} \big|\f(v-U_k)-\E_\pi(\f^U(v))\big| |b_k(U_k)|\bigg) \,dv \\
      & \leq 2 \int_{v=-1}^{T+1} \sqrt{(n-1)\var_\pi(\f(v-U))} \E\bigg(\sum_{k=1}^{n-1} \big|\f(v-U_k)-\E_\pi(\f(v-U))\big| |b_k(U_k)|\bigg) \,dv \\
      & \leq 2 \sqrt{(n-1)\frac{\|\f\|_2^2}{T}} \E\left(\sum_{k=1}^{n-1} \int_{v=-1}^{T+1} \big|\f(v-U_k)-\E_\pi(\f(v-U))\big| |b_k(U_k)| \,dv\right) \\
      & \leq 2 \sqrt{\frac{n-1}{T}} \|\f\|_2 \E\left(2\|\f\|_1 \sum_{k=1}^{n-1} |b_k(U_k)| \right) \\
      & \leq 4 \sqrt{\frac{n-1}{T}} \|\f\|_2 \|\f\|_1 \sqrt{n-1} \\
      & \leq 4 \frac{n-1}{\sqrt{T}} \|\f\|_1 \|\f\|_2,
      \end{align*}
      by using Lemma 6.1 of \cite{San12}. Then,
      \[D \leq C \frac{2^{-j/2}n}{\sqrt{T}},\]
      with $C$ an absolute positive constant;
  \item $\displaystyle B^2 = \sup_{u} \left(\sum_{k=1}^{n-1} \E(g^2(u,U_k))\right)$, with
      \begin{align*}
      & \E(g^2(u,U_k)) \\
      & = 4 \E\left[\left(\int_{v=-1}^{T+1} [\f(v-u)-\E_\pi(\f(v-U))][\f(v-U_k)-\E_\pi(\f(v-U))] \,dv\right)^2\right] \\
      & \leq 4 \E\left[\int_{v=-1}^{T+1} \big[\f^u(v)-\E_\pi(\f^U(v))\big]^2 \big|\f^{U_k}(v)-\E_\pi(\f^U(v))\big| \,dv \int_{v=-1}^{T+1} \big|\f^{U_k}(v)-\E_\pi(\f^U(v))\big| \,dv\right] \\
      & \leq 8 \E\left[\int_{v=-1}^{T+1} \big[\f(v-u)-\E_\pi(\f(v-U))\big]^2 \big|\f(v-U_k)-\E_\pi(\f(v-U))\big| \,dv\right] \|\f\|_1 \\
      & \leq \frac{16}{T} \int_{v=-1}^{T+1} \big[\f(v-u)-\E_\pi(\f(v-U))\big]^2 \,dv \|\f\|_1^2 \\
      & \leq \frac{64}{T} \|\f\|_1^2 \|\f\|_2^2,
      \end{align*}
      by using Lemma 6.3 of \cite{San12}. Hence,
      \[B^2 \leq C \frac{2^{-j}n}{T},\]
      with $C$ an absolute positive constant.
\end{itemize}
Finally, we obtain for all $\omega>0$, with probability larger than $1 - 2 \times 2.77 e^{-\omega}$,
\begin{equation}
|W_2| \leq C\left\{\frac{2^{-j/2}n}{\sqrt{T}}\sqrt{\omega} + \frac{2^{-j}n}{T}\sqrt{\omega} + \frac{2^{-j/2}n}{\sqrt{T}}\omega + \frac{2^{-j/2}\sqrt{n}}{\sqrt{T}}\omega^{3/2} + \omega^2\right\},  \label{pr-W2}
\end{equation}
with $C$ an absolute positive constant.

Thus, by inequalities (\ref{pr-majVS}), (\ref{pr-W0}), (\ref{pr-W1}) and (\ref{pr-W2}), for all $\omega>0$, with probability larger than $1 - 2 \times 2.77 e^{-\omega}$,
\begin{equation}
V_S \leq \frac{C(\omega)}{T} \left\{n + \frac{2^{-j}n^2}{T} + \frac{2^{-j/2}n}{\sqrt{T}}\right\}.  \label{pr-VS}
\end{equation}

\medskip

Then it remains to compute $B_S$. We recall that
\begin{align*}
B_S & = \sup_{v \in [-1;T+1]} \left|\sum_{i=1}^{n} \left[\f(v-U_i) - \frac{n-1}{n} \E_\pi(\f(v-U))\right]\right| \\
& \leq \tilde{B}_S + \frac{\|\f\|_1}{T},
\end{align*}
with $\displaystyle \tilde{B}_S = \sup_{v \in [-1;T+1]} \left|\sum_{i=1}^{n} \big[\f(v-U_i) - \E_\pi(\f(v-U))\big]\right|$.
Since the Haar basis is considered here, we can write for any $x\in\R$:
\[\f(x)=2^{j/2}\left(\1_{(2k+1)2^{-(j+1)} < x \leq (k+1)2^{-j}} - \1_{k2^{-j} \leq x \leq (2k+1)2^{-(j+1)}}\right),\]
with $\lambda=(j,k)$. Thus,
\[\tilde{B}_S \leq 2^{j/2} \left(\tilde{B}_S^1 + \tilde{B}_S^2\right),\]
where
\[\tilde{B}_S^1 = \sup_{v \in [-1;T+1]} \left|\sum_{i=1}^{n} \big[\1_{k2^{-j} \leq v-U_i \leq (2k+1)2^{-(j+1)}} - \E_\pi(\1_{k2^{-j} \leq v-U \leq (2k+1)2^{-(j+1)}})\big]\right|\]
and
\[\tilde{B}_S^2 = \sup_{v \in [-1;T+1]} \left|\sum_{i=1}^{n} \big[\1_{(2k+1)2^{-(j+1)} < v-U_i \leq (k+1)2^{-j}} - \E_\pi(\1_{(2k+1)2^{-(j+1)} < v-U \leq (k+1)2^{-j}})\big]\right|.\]

We observe that
\[\tilde{B}_S^1 \leq \sup_{B_v, v \in \R}\left|\sum_{i=1}^{n} \left[\1_{B_v}(U_i) - \E_\pi\big(\1_{B_v}(U)\big)\right]\right|,\]
where for any $v \in \R$, $B_v=[v-(2k+1)2^{-(j+1)};v-k2^{-j}]$.
We set $\mathcal{B}=\{B_v, v \in \R\}$ and for every integer $n$, $\displaystyle m_n(\mathcal{B}) = \sup_{A \subset \R, |A|=n} |\{A \cap B_v, v \in \R\}|$. It is easy to see that
\[m_n(\mathcal{B}) \leq 1 + \frac{n(n+1)}{2}\]
and so, the VC-dimension $\mathcal{V}$ of $\mathcal{B}$ defined by $\sup \{n \geq 0, m_n(\mathcal{B}) = 2^n\}$ is bounded by 2 (see Definition 6.2 of \cite{Mas07}).
Let us define $\sigma^2=\max\left\{2^{-(j+1)},K^2\mathcal{V}\big(1+\frac{j+1}{2}\ln{2}\big)/n\right\}$ with $K$ the absolute constant given by Lemma 6.4 of \cite{Mas07}. The quantity $\sigma^2$ satisfies in particular the two following assertions:
\[\forall B\in\mathcal{B}, \P_\pi[U \in B]\leq\sigma^2 \quad \mbox{and} \quad \sigma \geq K \sqrt{\mathcal{V}(1+\ln{(\sigma^{-1}\vee1)})/n}.\]
Indeed, if $\sigma^2=2^{-(j+1)}$, we have $K^2 \mathcal{V}(1+\ln{(\sigma^{-1}\vee1)})/n \leq K^2 \mathcal{V}(1+\ln{(2^{(j+1)/2})})/n \leq \sigma^2$, or else if $\sigma^2=K^2\mathcal{V}\big(1+\frac{j+1}{2}\ln{2}\big)/n$, we have $\sigma^{-1} \leq 2^{(j+1)/2}$ and so,
\[K^2 \mathcal{V}(1+\ln{(\sigma^{-1}\vee1)})/n \leq K^2 \mathcal{V}(1+\ln{(2^{(j+1)/2})})/n = \sigma^2.\]

By applying Lemma 6.4 of \cite{Mas07}, we obtain:
\begin{align*}
\E(\tilde{B}_S^1) & \leq \frac{K}{2} \sigma \sqrt{\mathcal{V}(1+|\ln{\sigma}|)} \\
& \leq \frac{K}{2} 2^{-(j+1)/2} \sqrt{\mathcal{V}\left(1+\frac{j+1}{2}\ln{2}\right)} + \frac{K}{2} K\mathcal{V}\left(1+\frac{j+1}{2}\ln{2}\right)/\sqrt{n} \\
& \leq C \left\{2^{-j/2}\sqrt{j} + \frac{j}{\sqrt{n}}\right\},
\end{align*}
with $C$ a positive absolute constant.
So, with a similar argument for $\tilde{B}_S^2$, we obtain for any $\lambda$ in $\Gamma$
\[\E(\tilde{B}_S) \leq C \left\{\sqrt{j} + \frac{j2^{j/2}}{\sqrt{n}}\right\}.\]
Consequently,
\[\E(B_S) \leq C \left\{\sqrt{j} + \frac{j2^{j/2}}{\sqrt{n}} + \frac{2^{-j/2}}{T}\right\},\]
with $C$ an absolute positive constant and from Markov's inequality, we have that for all $\omega>0$
\begin{equation}
\P\left(B_S > C(\omega) \left\{\sqrt{j} + \frac{j2^{j/2}}{\sqrt{n}} + \frac{2^{-j/2}}{T}\right\}\right) \leq e^{-\omega}.  \label{pr-BS}
\end{equation}

\medskip

Thus, by combining inequalities (\ref{def_f}), (\ref{pr-VS}) and (\ref{pr-BS}), we obtain for all $\omega>0$, with probability larger than $1 - (1+2 \times 2.77) e^{-\omega}$,
\begin{align*}
&f(U_1,\ldots,U_n;m) \\
& \leq \frac{C(\omega)}{n} \Bigg\{\sqrt{m \ln{(2/\alpha)} \left(\frac{n}{T} + \frac{2^{-j}n^2}{T^2} + \frac{2^{-j/2}n}{T^{3/2}}\right)} + \ln{(2/\alpha)} \left(\sqrt{j} + \frac{j2^{j/2}}{\sqrt{n}} + \frac{2^{-j/2}}{T}\right)\Bigg\}.
\end{align*}

Furthermore, $N_{[-1;T+1]}\sim\mathcal{P}((T+2)\mu_c+n\|h\|_1)$. Hence,
\[\E(N_{[-1;T+1]})\leq C(\mu_c,R_1) (n+T).\]
From Markov's inequality, we have that for all $\omega>0$
\[\P\left(N_{[-1;T+1]} > C(\omega,\mu_c,R_1) (n+T)\right) \leq e^{-\omega}.\]

Then, we choose $\omega$ such that this quantity $(2 \times 2.77 + 2) e^{-\omega}$ is equal to $\beta/2$.
So, with probability larger than $1 - \beta/2$,
\begin{align*}
& f(U_1,\ldots,U_n;m) \\
& \leq \frac{C(\beta,\mu_c,R_1)}{n} \Bigg\{\sqrt{\ln{(2/\alpha)}} \sqrt{(n+T) \left(\frac{n}{T} + \frac{2^{-j}n^2}{T^2} + \frac{2^{-j/2}n}{T^{3/2}}\right)} + \ln{(2/\alpha)} \left(\sqrt{j} + \frac{j2^{j/2}}{\sqrt{n}} + \frac{2^{-j/2}}{T}\right)\Bigg\} \\
& \leq \frac{C(\beta,\mu_c,R_1)}{n} \left\{\sqrt{\ln{(2/\alpha)}} \sqrt{n + \frac{n^2}{T} + \frac{2^{-j}n^3}{T^2}} + \ln{(2/\alpha)} \left(\sqrt{j} + \frac{j2^{j/2}}{\sqrt{n}} + \frac{2^{-j/2}}{T}\right)\right\} \\
& \leq C(\beta,\mu_c,R_1) \left\{\sqrt{\ln{(2/\alpha)}} \left(\frac{1}{\sqrt{n}} + \frac{1}{\sqrt{T}} + \frac{2^{-j/2}\sqrt{n}}{T}\right) + \ln{(2/\alpha)} \left(\frac{\sqrt{j}}{n} + \frac{j2^{j/2}}{n^{3/2}} + \frac{2^{-j/2}}{nT}\right)\right\}.
\end{align*}
Therefore by definition of $q^{\alpha}_{1-\beta/2}$,
\[q^{\alpha}_{1-\beta/2} \leq C(\beta,\mu_c,R_1) \left\{\sqrt{\ln{(2/\alpha)}} \left(\frac{1}{\sqrt{n}} + \frac{1}{\sqrt{T}} + \frac{2^{-j/2}\sqrt{n}}{T}\right) + \ln{(2/\alpha)} \left(\frac{\sqrt{j}}{n} + \frac{j2^{j/2}}{n^{3/2}} + \frac{2^{-j/2}}{nT}\right)\right\},\]
which concludes the proof of \ref{lem_quant}.

\bigskip

\noindent \textbf{Acknowledgments:}
The authors are very grateful to Patricia Reynaud-Bouret and Vincent Rivoirard for stimulating discussions and constructive comments and they would like to thank the anonymous Associate Editor and two referees for their helpful comments and valuable suggestions.
They wish also thank Magalie Fromont and her coauthors for their programs for the test of homogeneity.
The research of the authors is partly supported by the french Agence Nationale de la Recherche (ANR 2011 BS01 010 01 projet Calibration).

\bibliographystyle{plain}
\bibliography{Ref_PoissInterDetection}

\end{document}